\documentclass{template/cocv}
\usepackage[T1]{fontenc}
\usepackage[utf8]{inputenc}
\usepackage[english]{babel}

\usepackage{amssymb}
\usepackage{amsthm}
\usepackage{amsmath}

\usepackage[backend=bibtex,%
			style=numeric,%
			autolang=hyphen,%
			isbn=false,%
			date=year,%
			eprint=false,%
			]{biblatex}

\usepackage{xcolor}
\usepackage{enumitem}
\usepackage{cancel}
\usepackage{dirtytalk}
\usepackage{csquotes}
\usepackage{float,tabularx}

\usepackage{hyperref} 
\usepackage[nameinlink,capitalise,noabbrev]{cleveref}

\usepackage{soul}

\hypersetup{
	breaklinks=true,
    colorlinks,
    linkcolor={blue!80!black},
    citecolor={blue!80!black},
    urlcolor={blue!80!black},
    pdftitle={Trivialisable control-affine systems revisited},
    pdfauthor={Timothée Schmoderer, Witold Respondek},
    pdfsubject={Article on trivial control-systems for ESAIM: Control, Optimisation and Calculus of Variations},
    pdfkeywords={Control-affine system, Feedback equivalence, Trivial control systems, Control curvature, Normal forms, Infinitesimal symmetries},
}

\AtEveryBibitem{\clearfield{pages}}    
\AtEveryBibitem{\clearfield{edition}}  
\AtEveryBibitem{\clearfield{place}}    
\AtEveryBibitem{\clearfield{address}}  
\AtEveryBibitem{\clearfield{location}} 
\AtEveryCitekey{\clearlist{location}}
\AtEveryBibitem{\clearlist{location}}
\renewbibmacro{in:}{}				   

\setstcolor{red}

\newcommand{\NN}{\ensuremath{\mathbb{N}}}

\newcommand{\RR}{\ensuremath{\mathbb{R}}}

\newcommand{\AAA}{\ensuremath{\mathcal{A}}}

\newcommand{\FFF}{\ensuremath{\mathcal{F}}}
\newcommand{\GGG}{\ensuremath{\mathcal{G}}}

\newcommand{\MMM}{\ensuremath{\mathcal{M}}}

\newcommand{\XXX}{\ensuremath{\mathcal{X}}}

\newcommand{\AAAA}{\ensuremath{\mathfrak{A}}}
\newcommand{\IIII}{\ensuremath{\mathfrak{I}}}
\newcommand{\LLLL}{\ensuremath{\mathfrak{L}}}

\crefformat{equation}{(#2#1#3)}
\crefmultiformat{equation}{equations (#2#1#3)}{ and ~(#2#1#3)}{, (#2#1#3)}{ and~(#2#1#3)}
\crefrangelabelformat{equation}{equations (#3#1#4) to~(#5#2#6)}
\crefformat{enumi}{#2#1#3}
\crefmultiformat{enumi}{#2#1#3}{ and ~#2#1#3}{, #2#1#3}{ and~#2#1#3}
\crefname{problem}{problem}{problems}
\Crefname{problem}{Problem}{Problems}

\crefformat{figure}{#2figure~#1#3}
\Crefformat{figure}{#2Figure~#1#3}

\newcommand{\lccref}[1]{\hyperref[{#1}]{\lcnamecref{#1}~\labelcref{#1}}} 

\newcommand\lb[2]{\left[#1,#2\right]}              
\renewcommand\vec[1]{\frac{\partial}{\partial #1}} 
\newcommand\distrib[1]{\mathrm{span}\left\{#1\right\}} 
\newcommand\ad[2]{\mathrm{ad}_{#1}#2}                
\newcommand\adk[3]{\mathrm{ad}_{#1}^{#2}#3}          
\newcommand\dL[2]{\mathrm{L}_{#1}\left(#2\right)}  
\newcommand\dLk[3]{\mathrm{L}_{#1}^{#2}\left(#3\right)\,} 
\newcommand\vectR[1]{\mathrm{vect}_{\RR}\left\{#1\right\}}       
\newcommand\sgn[1]{\mathrm{sgn}\left(#1\right)}      

\newcommand\rk[1]{\mathrm{rk}\,#1}        

\newcommand\diff{\mathrm{d}} 

\makeatletter
\newcommand*\owedge{\mathpalette\@owedge\relax}
\newcommand*\@owedge[1]{%
  \mathbin{%
    \ooalign{%
      $#1\m@th\bigcirc$\cr
      \hidewidth$#1\m@th\wedge$\hidewidth\cr
    }%
  }%
}
\makeatother

\newtheorem{theorem}{Theorem}
\newtheorem{definition}{Definition}
\newtheorem{proposition}{Proposition}
\newtheorem{lemma}{Lemma}

\newtheorem{remark}{Remark}

\newcommand{\myparagraph}[1]{\paragraph{\textbf{#1}}}
\begin{document}
\title{Trivialisable control-affine systems revisited}
\runningtitle{Trivialisable control systems}
\author{Timothée Schmoderer}\address{Laboratoire de Mathématiques de l'INSA UR 3226 - FR CNRS 3335, INSA Rouen Normandie, Avenue de l'Université, 76800 St Etienne du Rouvray, France;
\email{\href{mailto:timothee.schmoderer@insa-rouen.fr}{timothee.schmoderer@insa-rouen.fr} \& \href{mailto:witold.respondek@insa-rouen.fr}{witold.respondek@insa-rouen.fr}}}
\author{Witold Respondek}\sameaddress{1}\secondaddress{Institute of Automatic Control, {\L}\'{o}d\'{z} University of Technology, Poland}
\date{\today}
\begin{abstract}
The purpose of this paper is to explore the concept of trivial control systems, namely systems whose dynamics depends on the controls only. Trivial systems have been introduced and studied by Serres in the the context of control-nonlinear systems on the plane with a scalar control. In our work, we begin by proposing an extension of the notion of triviality to control-affine systems with arbitrary number of states and controls. Next, our first result concerns two novel characterisations of trivial control-affine systems, one of them is based on the study of infinitesimal symmetries and is thus geometric. Second, we derive a normal form of trivial control-affine systems whose Lie algebra of infinitesimal symmetries possesses a transitive almost abelian Lie subalgebra. Third, we study and propose a characterisation of trivial control-affine systems on $3$-dimensional manifolds with scalar control. In particular, we give a novel proof of the previous characterisation obtained by Serres. Our characterisation is based on the properties of two functional feedback invariants: the curvature (introduced by Agrachev) and the centro-affine curvature (used by Wilkens). Finally, we give several normal forms of control-affine systems, for which the curvature and the centro-affine curvature have special properties.
\end{abstract}
\begin{resume}
L'objectif de cet article est d'explorer le concept de système de contrôles trivial,c'est-à-dire un système dont la dynamique dépend uniquement des commandes. Les systèmes triviaux ont été introduits et étudiés par Serres dans le contexte des systèmes bidimensionnels et non-linéaires par rapport à un contrôle scalaire. Dans notre travail, nous commençons par proposer une extension de la notion de trivialité aux systèmes affines par rapport aux contrôles et avec un nombre arbitraire d'états et de contrôles. Ensuite, notre premier résultat concerne deux nouvelles caractérisations des systèmes affines triviaux, l'une d'entre elles est basée sur l'étude des symétries infinitésimales et est donc géométrique. Deuxièmement, nous donnons une forme normale pour les systèmes affines triviaux dont l'algèbre de Lie des symétries infinitésimales possède une sous-algèbre de Lie presque abélienne. Troisièmement, nous étudions et proposons une caractérisation des systèmes affines triviaux sur des variétés à $3$ dimensions et avec un contrôle scalaire. En particulier, nous donnons une nouvelle preuve de la caractérisation précédemment obtenue par Serres. Notre caractérisation est basée sur les propriétés de deux invariants fonctionnels du bouclage : la courbure (introduite par Agrachev) et la courbure centro-affine (utilisée par Wilkens). Enfin, nous donnons plusieurs formes normales de systèmes affine, pour lesquelles la courbure et la courbure centro-affine ont des propriétés particulières.
\end{resume}
\maketitle
\section{Introduction}
In this paper, we consider control-affine systems $\Sigma$ of the form
\begin{align}\label{eq:control_affine_system_nm}
    \Sigma\,:\,\dot{\xi}=f(\xi)+\sum_{i=1}^mg_i(\xi)u_i,\quad u_i\in\RR,
\end{align}
\noindent
where the state $\xi$ belongs to a smooth $n$-dimensional manifold $\MMM$ (or an open subset of $\RR^n$, since most of our results are local), and $f$ and $g$ are smooth vector fields on $\MMM$, i.e. smooth sections of the tangent bundle $T\MMM$. {Throughout the paper, the word "smooth" will always mean $C^{\infty}$-smooth, and all objects (manifolds, vector fields, differential forms, functions) are assumed to be smooth.} We denote a control-affine system by the pair $\Sigma=(f,g)$, where $g=(g_1,\ldots,g_m)$. To any control-affine system $\Sigma=(f,g)$ we attach two distributions: 
\begin{align}\label{eq:def_distributions_g_g1}
    \GGG=\distrib{g_1,\ldots,g_m}\quad\textrm{and}\quad\GGG^1=\GGG+\lb{f}{\GGG}=\distrib{g_1,\ldots,g_m,\lb{f}{g_1},\ldots,\lb{f}{g_m}}.
\end{align}
\noindent
%
We call two control-affine systems $\Sigma=(f,g)$ and $\tilde{\Sigma}=(\tilde{f},\tilde{g})$ \emph{feedback equivalent}, if there exists a diffeomorphism $\phi :\MMM\rightarrow\tilde{\MMM}$ and smooth functions $\alpha:\MMM\rightarrow\RR^m$ and $\beta:\MMM\rightarrow GL_m(\RR)$ such that 
\begin{align*}
    \tilde{f}=\phi_*\left(f+\sum_{i=1}^m g_i\alpha_i\right)\quad\textrm{and}\quad\tilde{g}_i=\phi_* \left(\sum_{j=1}^m g_j\beta_i^j\right),
\end{align*}
\noindent
where $\phi_*$ denotes the tangent map of $\phi$. If $\phi$ is defined locally around $\xi_0$ and $\tilde{\xi}_0=\phi(\xi_0)$, then we say that $\Sigma$ and $\tilde{\Sigma}$ are locally feedback equivalent at $\xi_0$ and $\tilde{\xi}_0$, respectively. 
Feedback equivalence of control-affine systems means equivalence of the affine distributions $\AAA=f+\GGG$ and $\tilde{\AAA}=\tilde{f}+\tilde{\GGG}$ attached to $\Sigma$ and $\tilde{\Sigma}$, respectively. \\

\noindent
In the thesis \cite{serres2006Geometryfeedbackclassification}, Serres proposed the notion of a trivial system of the form 
\begin{align*}
(\mathcal{T})\,:\,\dot{x}=F(w),\quad x\in\XXX,\quad w\in\RR^m,
\end{align*}
\noindent
where $w$ is the control that enters nonlinearly. The dynamics $F(w)$ of a trivial system does not depend on the state variables $x$ and thus depends on control variables $w$ only. Actually, $(\mathcal{T})$ is called flat in \cite{serres2006Geometryfeedbackclassification} but that name can be misleading because{, first, there is a well established notion of flat {control} systems \cite{fliess1992DifferentiallyFlatNonlinear} and, second,} the class of trivial control systems does not coincide with control systems of zero-curvature \cite{agrachev1997Feedbackinvariantoptimalcontrol}{, which thus can be considered as geometrically flat}, as we will discuss in \cref{sec:trivial_systems_3d}. {For those reasons, following \cite{serres2007curvaturefeedbackclassification}, we call $(\mathcal{T})$ a trivial system and we say that a} general control-nonlinear system $\dot{x}=F(x,w)$ is  trivialisable if it is equivalent, via a feedback of the form $\tilde{x}=\phi(x)$, $\tilde{w}=\psi(x,w)$, to a trivial system $(\mathcal{T})$, where $(\phi,\psi):\XXX\times\RR^m\rightarrow\tilde{\XXX}\times\RR^m$ is a diffeomorphism. {Inspired by the above considerations,} we adapt the concept of triviality to control-affine systems as follows.
\begin{definition}[Trivial {control-affine} systems]\label{def:trivial_system}
We say that a control-affine system $\Sigma=(f,g)$ is trivialisable if it is feedback equivalent {to a trivial system of the form:}
\begin{align*}
    (T)\,:\,\left\{
    \begin{array}{rl}
        \dot{x} &= F(w) \\
        \dot{w} &= u 
    \end{array}
    \right., \quad (x,w)\in\MMM=\XXX\times\RR^m,\quad u\in\RR^m,
\end{align*}
\noindent
{whose} {$\dot{x}$-}dynamics depend on the controlled {$w$}-variables {only}. 
\end{definition}
{The notions of trivial and trivialisable general control-nonlinear {versus} control-affine systems are two sides of the same coin. Indeed, two control-nonlinear systems $\dot{x}=F(x,w)$ and $\dot{\tilde{x}}=\tilde{F}(\tilde{x},\tilde{w})$ are feedback equivalent if and only if their control-affine extensions $\dot{x}=F(x,w),\,\dot{w}=u$ and $\dot{\tilde{x}}=\tilde{F}(\tilde{x},\tilde{w}),\,\dot{\tilde{w}}=\tilde{u}$ are equivalent via control-affine feedback transformations, see \cite[equation 3.6]{jakubczyk1990EquivalenceInvariantsNonlinear}. Therefore a control-nonlinear system $\dot{x}=F(x,w)$ is trivialisable if and only if $\dot{x}=F(x,w),\,\dot{w}=u$ is trivialisable in the sense of \cref{def:trivial_system} and the latter class is the object of our studies in this paper. }  \\

Trivial control systems are interesting to study because they model trajectories of dynamical systems under a nonholonomic constraint that does not depend on the point. Indeed, under the additional {regularity} assumption that $\rk\frac{\partial F}{\partial w}(w)=m$, equivalently, the distribution $\GGG^1$ of $(T)$ satisfies $\rk\GGG^1=2m$, there exist local coordinates $x=(z,y)$, with $\dim z=n-2m$ and $y=(y_1,\ldots,y_m)$ such that the equations of $(T)$ can be rewritten 
\begin{align*}
    \left\{\begin{array}{rl}
        \dot{z} &= \textsf{f}(w)  \\
        \dot{y} &= w \\
        \dot{w} &= u
    \end{array}\right.
\end{align*}
\noindent
and {we conclude} that the trajectories of $(T)$ satisfy {the nonholonomic constraints} $\dot{z}=\textsf{f}(\dot{y})${, whose shape is independent of the point $x=(z,y)$}. Denoting by $\XXX$ the (locally defined) quotient manifold $\MMM/\GGG$, we see that a trajectory $x(t)\in\XXX$ satisfies the nonholonomic constraint $\dot{z}=\textsf{f}(\dot{y})$ if and only if there exists a smooth control $u(t)$ such that $(x(t),w(t))$ is a trajectory of $(T)$. Connections between equations on the tangent bundle and control systems are explored in \cite{schmoderer2018Studycontrolsystems,schmoderer2021Conicnonholonomicconstraints}. Examples of trivial systems can be found in the literature; e.g. in \cite{schmoderer2021Conicnonholonomicconstraints} we characterise trivial elliptic, hyperbolic, and parabolic control systems, Dubin's car \cite{dubins1957CurvesMinimalLength} is a very simple model of system that is trivial, and, finally, trivial control-nonlinear system on surfaces (i.e. $n=2$) and with scalar control have been studied, characterised (and normal forms in particular cases have been given) in \cite{serres2006Geometryfeedbackclassification,serres2007curvaturefeedbackclassification,serres2009Controlsystemszero}.

\subsection{Outline of the paper}
In the next subsection, we develop the main notions of differential geometry and of control theory that we will need in the rest of the paper. Next, in \cref{sec:trivial_systems_general}, we study trivial control-affine systems on manifolds of arbitrary dimension and with an arbitrary number of controls. We propose two novel characterisations of trivial systems, one of them is based on the Lie algebra of infinitesimal symmetries. Moreover, using our characterisation of trivial systems via symmetries, we will give a normal form of trivial systems whose Lie algebra of infinitesimal symmetries possesses a transitive almost abelian Lie subalgebra. Afterwards, in \cref{sec:trivial_systems_3d}, we will be interested in revisiting the characterisation of trivial systems discovered by {Serres} \cite{serres2006Geometryfeedbackclassification} in the context of control-nonlinear systems on surfaces. We propose a characterisation of trivial control-affine systems on $3$-dimensional manifolds with scalar control. Our characterisation exhibits a discrete invariant, and two fundamental functional invariants: the control curvature introduced by {Agrachev}  \cite{agrachev1997Feedbackinvariantoptimalcontrol,agrachev1998FeedbackInvariantOptimalControl}, and the centro-affine curvature {used by Wilkens \cite{wilkens1998Centroaffinegeometryplane}}. Both functional invariants can be computed for any control-affine system. We will provide another proof and interpretation of Serres results. Finally, in \cref{sec:normal_forms} we discuss several normal forms {(some new and some existing in the litterature)} of control-affine systems, for which the control curvature and the centro-affine curvature have special properties.  
\subsection{Preliminaries}
In this subsection, we recall the main definitions and notions of differential geometry and of control theory that we need in the paper. The main notations that we use are summarised in \cref{table:main_notations}.\\
\myparagraph{Differential Geometry.} For a manifold $\MMM$ we will denote by $T\MMM$ and $T^*\MMM$ the tangent and cotangent bundle, respectively. The space of all smooth vector fields (smooth sections of $T\MMM$) will be denoted $V^{\infty}(\MMM)$ and the space of all smooth differential $p$-forms by $\Lambda^p(\MMM)$, except for smooth functions ($0$-forms) whose space is denoted $C^{\infty}(\MMM)$. For a diffeomorphism $\phi\,:\,\MMM\rightarrow\tilde{\MMM}$, a vector field $f\in V^{\infty}(\MMM)$, and a differential $p$-form $\omega\in\Lambda^p(\tilde{\MMM})$, we denote by $\phi_*f \in V^{\infty}(\tilde{\MMM})$ the push-forward of $f$, and  by $\phi^*\omega\in\Lambda^p(\MMM)$ the pull-back of $\omega$. The (local) flow of a vector field $f\in V^{\infty}(\MMM)$ is denoted by $\gamma_t^f$ (for any $t$ for which it is defined). The Lie derivative of a differential $p$-form $\omega$ along a vector field $f$ will be denoted by $\dL{f}{\omega}$. In particular, for a function $\lambda\in C^{\infty}(\MMM)$ and its differential $\diff\lambda$ (an exact $1$-form) we have 
\begin{align*}
\dL{f}{\lambda}=\langle\diff\lambda, f\rangle\quad\textrm{and}\quad\dL{f}{\diff\lambda}=\diff\dL{f}{\lambda}.
\end{align*}
\noindent
For any smooth functions $\alpha$, $\lambda$, and $\mu$, the Lie derivative possesses the following properties: $\dL{\alpha f}{\lambda}=\alpha\dL{f}{\lambda}$, and $\dL{f}{\lambda\mu}=\dL{f}{\lambda}\mu+\lambda\dL{f}{\mu}$. Iterative Lie derivatives are defined by $\dLk{f}{k}{\lambda}=\dL{f}{\dLk{f}{k-1}{\lambda}}$, for any $k\geq2$. For any two vector fields $f,g\in V^{\infty}(\MMM)$, we define their Lie bracket as a new vector field, denoted $\lb{f}{g}\in V^{\infty}(\MMM)$, such that for any smooth function $\lambda$ we have 
\begin{align*}
\dL{\lb{f}{g}}{\lambda}=\dL{f}{\dL{g}{\lambda}}-\dL{g}{\dL{f}{\lambda}}.
\end{align*}
\noindent
The Lie bracket possesses the following properties: it is bilinear over $\RR$, it is skew-commutative, i.e. $\lb{f}{g}=-\lb{g}{f}$, and it satisfies the Jacobi identity: 
\begin{align*}
\lb{f}{\lb{g}{h}}+\lb{g}{\lb{h}{f}}+\lb{h}{\lb{f}{g}} &= 0,\quad \forall\,f,g,h\in V^{\infty}(\MMM).
\end{align*}
\noindent
Moreover, for any smooth function $\alpha$, and any vector fields $f$, $g$, and $h$, we have 
\begin{align*}
\lb{f}{\alpha g+h}=\alpha\lb{f}{g}+\dL{f}{\alpha}g+\lb{f}{h}.
\end{align*}
\noindent
Two vector fields $f$ and $g$ satisfying $\lb{f}{g}=0$ are said to be commuting; since under diffeomorphisms $\phi\,:\,\MMM\rightarrow\tilde{\MMM}$ the Lie bracket is transformed by $\lb{\phi_*f}{\phi_*g}=\phi_*\lb{f}{g}$, the commutativity property does not depend on coordinates. The celebrated \emph{Flow-box} theorem (also called the "{Straightening-out theorem}" or the "{Local linearisation lemma}") asserts that on a given $n$-dimensional manifold $\MMM$ there exists a local coordinate system $(x_1,\ldots,x_n)$ such that $f=\vec{x_1}$ in a neighbourhood of {any} point $p$ where $f(p)\neq0$. This can simultaneously be done for a family of (locally) independent vector fields $(f_1,\ldots,f_m)$ if and only if they are mutually commuting. {We set $\adk{f}{0}{g}=g$, $\ad{f}{g}=\lb{f}{g}$, and } the iterated Lie bracket is denoted by $\adk{f}{k}{g}=\lb{f}{\adk{f}{k-1}{g}}$ for $k\geq 1$; see \cite[chapter 1]{isidori1995NonlinearControlSystems} for a detailed introduction and proofs of the above properties.

\myparagraph{Infinitesimal symmetries.}
We briefly introduce the notion of symmetries of control-affine system{s} (see \cite{respondek2002Nonlinearizablesingleinputcontrol,grizzle1985structurenonlinearcontrol} for a detailed introduction). For a control-affine system $\Sigma=(f,g)${, see \cref{eq:control_affine_system_nm}}, with state $\MMM$ a smooth $n$-dimensional manifold, we define the field of admissible velocities $\AAA$ as
\begin{align*}
    \AAA(\xi)=\{f(\xi)+\sum_{i=1}^m g_i(\xi)u_i\, :\, u_i\in\RR\}\subset T_{\xi}\MMM .
\end{align*}
\noindent
We say that a diffeomorphism $\phi:\MMM\rightarrow\MMM$ is a \emph{symmetry} of $\Sigma$ if it preserves the field of affine $m$-planes $\AAA$ (equivalently, the affine distribution $\AAA=f+\GGG$), that is, $\phi_*\AAA=\AAA$. 
We say that a vector field $v$ on $\MMM$ is an infinitesimal symmetry of $\Sigma=(f,g)$ if the (local) flow $\gamma^v_t$ of $v$ is a local symmetry, for any $t$ for which it exists, that is, $(\gamma_t^v)_*\AAA=\AAA$. Consider the system $\Sigma=(f,g)$ and recall that $\GGG$ is the distribution spanned by the vector fields $g_1,\ldots,g_m$. We have the following characterisation of infinitesimal symmetries. 
\begin{proposition}\label{prop:charac_lie_algebra_infitesimal_symmetries}
A vector field $v$ is an infinitesimal symmetry of the control-affine system $\Sigma=(f,g)$ if and only~if 
\begin{align*}
\lb{v}{g}=0\!\!\mod\GGG\quad\textrm{and}\quad \lb{v}{f}=0\!\!\mod\GGG.
\end{align*}
\end{proposition}
\noindent
By the Jacobi identity, it is easy to see that if $v_1$ and $v_2$ are infinitesimal symmetries, then so is $\lb{v_1}{v_2}$, hence the set of all infinitesimal symmetries forms a real Lie algebra. Notice that the Lie algebra of infinitesimal symmetries is attached to the affine distribution $\AAA=f+\GGG$ and not to {a particular {pair} $(f,g)=(f,g_1,\ldots,g_m)$. Different pairs $(f,g)$ related via feedback transformations $(\alpha,\beta)$ define the same $\AAA$ and thus {have} the same Lie algebra of infinitesimal symmetries which, therefore, is a feedback invariant object attached to $\Sigma$.}

\begin{table}[H]
\centering
\setlength{\tabcolsep}{1em} 
\renewcommand{\arraystretch}{1.2}
\begin{tabularx}{\textwidth}{lX}
{$\MMM$, $T\MMM$, $\xi=(x,w)$} &  {Smooth $n$-dimensional manifold, its tangent bundle, and its local coordinates with $\dim w=m$.} \\
{$\phi$, $\phi_*$, $\phi^*$}  & {A diffeomorphism, its tangent map, its cotangent map.}  \\
$\Sigma=(f,g)$ & A control-affine system; see \cref{eq:control_affine_system_nm}. \\
$\GGG$ and $\GGG^1$ & Distribution spanned by the vector fields $g_1,\ldots, g_m$ and {the} distribution spanned by the vector fields $g_1,\ldots,g_m$ and $\lb{f}{g_1},\ldots,\lb{f}{g_m}$; see \cref{eq:def_distributions_g_g1}. \\
$(T)$ & Trivial control-affine system; see \cref{def:trivial_system}. \\
%
$\LLLL$, $\AAAA$, $\IIII$ & A real Lie (sub)algebra, a subalgebra, an ideal. \\ 
$\Sigma_{\lambda}$, $\Sigma_{\lambda}^{0,k}$ & Normal forms of trivial systems having an almost abelian subalgebra of infinitesimal symmetries; see \cref{thm:normal_form_trivial_almost_abelian,prop:normal_form_trivial_system_non_generic}. \\ 
$\Sigma_s=(f_s,g)$ & Control affine system given by a semi-canonical pair; see \cref{def:pair_control_affine_system}. \\ 
$\Sigma_c=(f_c,g_c)$ & Control affine system given by the canonical pair; see \cref{def:pair_control_affine_system}.\\ 
$(k_1,k_2,k_3)$ and $(\lambda_1,\lambda_2,\lambda_3)$ & Structure functions attached to any control-affine system on a $3$-dimensional manifold with scalar control; see \cref{eq:structure_functions_definition}. \\ 
$(\varepsilon,\kappa,\nu)$ & Feedback invariants of control-affine systems; defined for the canonical pair by \cref{eq:canonical_structure_constants} and expressed for any control-affine system by \cref{eq:formula_curavtures_struct_functions}.
\end{tabularx}%
\caption{Main notations for the paper}\label{table:main_notations}
\end{table}%

\section{Trivial control-affine systems}\label{sec:trivial_systems_general}

In this section, we first propose two new characterisations of trivialisable control-affine systems (with the state-space of arbitrary dimension and with an arbitrary number of controls); see \cref{thm:charac_trivial_sys} below. Second, we give a normal form of trivial systems whose Lie algebra of infinitesimal symmetries possesses a transitive almost abelian Lie subalgebra; see \cref{thm:normal_form_trivial_almost_abelian,prop:normal_form_trivial_system_non_generic} of this section.

\subsection{Characterisations of trivial systems}

The following theorem gives two characterisations of trivialisable systems. The first one is technical and shows that triviality is a property that depends on the coordinates (like being a linear control system depends on the choice of coordinates), and the second one is based on infinitesimal symmetries and is thus geometric. Recall that to a control-affine system $\Sigma=(f,g)$ we attach two distributions $\GGG=\distrib{g_1,\ldots,g_m}$ and $\GGG^1=\GGG+\lb{f}{\GGG}$, see \cref{eq:def_distributions_g_g1}.

\begin{theorem}[Two characterisations of trivialisable systems]\label{thm:charac_trivial_sys}\leavevmode
Consider a control-affine system $\Sigma=(f,g)$ with state on a $n$-dimensional manifold $\MMM$ and with $m\geq1$ controls. The following assertions hold locally around $\xi_0$:
\begin{enumerate}[label=(\roman*),ref=\textit{(\roman*)}]
\item \label{thm:charac_trivial_sys:1} {Suppose that $\rk\GGG^1=m+k$ is constant. The system $\Sigma$ is locally trivialisable if and only if $\Sigma$ is locally feedback equivalent to 
    \begin{align*}
        \Sigma_T\,:\,\left\{
        \begin{array}{rl}
            \dot{x}_i &= h_i(x,w) \\
            \dot{w} &= u
        \end{array}
        \right.,\quad \textrm{for}\quad 1\leq i\leq n-m,
    \end{align*}
    \noindent
    where the smooth scalar functions $h_1,\ldots,h_{n-m}$ satisfy
    \begin{align}\label{thm:charac_trivial_sys:eq:condition_h}
        \rk\distrib{\diff h_1,\ldots,\diff h_{n-m}}=k.
    \end{align}}
    \item\label{thm:charac_trivial_sys:2} $\Sigma$ is, locally around $\xi_0$, trivialisable if and only if the distribution $\GGG$ is involutive and of constant rank $m$ and, additionally, the Lie algebra of infinitesimal symmetries of $\Sigma$ possesses an abelian subalgebra $\AAAA$ such that $\AAAA(\xi_0)\oplus\GGG(\xi_0)=T_{\xi_0}\MMM$.
\end{enumerate}
\end{theorem}
Observe that the assumption on the rank of the distribution $\GGG^1$ in statement \cref{thm:charac_trivial_sys:1} implies that the dimension $n$ of the manifold $\MMM$ is greater than or equal to $m+k$. If $n=m+k$, then the trivialisation $(T)$ of $\Sigma_T$ (and thus of $\Sigma$) can be taken (for suitable $w$ and $u$) as $\dot{x}_i=w_i$, $1\leq i\leq k$, $\dot{w}_j=u_j$, $1\leq j\leq m$. On the other hand, if $n>m+k$, then $\dot{x}_i=w_i$, $1\leq i\leq k$, $\dot{x}_i=F_i(w_1,\ldots,w_k)$, $k+1\leq i\leq m+k$, and $\dot{w}_j=u_j$, $1\leq j\leq m$. Notice that $k\leq m$, so if $n > 2m$, then there are always {nonlinear} equations $\dot{x}_i=F_i(w_1,\ldots,w_k)$. In item \cref{thm:charac_trivial_sys:2}, there are no particular relations between the dimension of the state space and on the number of controls (other than the obvious $n\geq m$).
\begin{remark}
For the system $\Sigma_T$, define $H=(h_1,\ldots,h_{n-m})^t$. Then, under the assumption that $\rk\GGG^1$ is constant, condition \cref{thm:charac_trivial_sys:eq:condition_h} can be equivalently reformulated as 
\begin{align*}
    \rk\frac{\partial H}{\partial w}(x,w) = \rk\frac{\partial H}{\partial (x,w)}(x,w),
\end{align*}
\noindent
in a neighbourhood of $(x_0,w_0)$.
\end{remark}
\begin{proof}\leavevmode
\begin{enumerate}[label=\textit{(\roman*)}]
	\item Suppose that $\Sigma$ is locally trivialisable, {i.e. by \cref{def:trivial_system},} $\Sigma$ is locally feedback equivalent to $(T)$, which is of the form of $\Sigma_T$ with $h_i(x,w)=F_i(w)$, for $1\leq i\leq n-m$, and we now show that those functions satisfy \cref{thm:charac_trivial_sys:eq:condition_h}. On one hand, the condition $\rk\GGG^1=m+k$ implies that the Jacobian matrix $\frac{\partial F}{\partial w}$ {is of constant rank $k$}, where $F=(F_1,\ldots,F_{n-m})^T$. On the other hand, we obtain that ${\diff h_i}=\diff F_i=\sum_{j=1}^{m} \frac{\partial F_i}{\partial w_j}\diff w_j$. Hence the rank of $\distrib{\diff h_1,\ldots,\diff h_{n-m}}$ is the same as {that} of $\frac{\partial F}{\partial w}$ and the conclusion follows. {Conversely, assume that $\Sigma$ is feedback equivalent to $\Sigma_T$. Using the assumption $\rk\GGG^1=m+k$ we can reorder the $x$-coordinates such that $h=(\hat{h}_1,\ldots,\hat{h}_k,\tilde{h}_{k+1},\ldots,\tilde{h}_{n-m})$, where $\rk\frac{\partial \hat{h}}{\partial w}=k$. {We set $\hat{w}_i=\hat{h}_i(x,w)$, for $1\leq i\leq k$, completed by $\hat{w}_{k+1},\ldots,\hat{w}_m$ (chosen among the $w_i$'s) in such a way that $\hat{w}_1,\ldots,\hat{w}_m$ form a local coordinate system.} {W}e conclude, by condition \cref{thm:charac_trivial_sys:eq:condition_h}, that the functions $\tilde{h}_{k+1},\ldots,\tilde{h}_{n-m}$ depend on {the} variables $\hat{w}_1,\ldots,\hat{w}_k$ {only}. {Using a feedback transformation that yields $\dot{\hat{w}}_i=\hat{u}_i$, $1\leq i\leq k$, we conclude that} $\Sigma_T$ is{, indeed,} a trivial system in coordinates $(x,\hat{w})$.}
    \item Suppose that $\Sigma=(f,g)$ is {locally} trivialisable, then for $(T)$ we have $\GGG=\distrib{\vec{w_1},\ldots,\vec{w_m}}$, which clearly is involutive and of constant rank $m$. Moreover, the vector fields $v_i=\vec{x_i}$, for $1\leq i\leq n-m$, are commuting symmetries of $(T)$ {that span the abelian Lie algebra $\AAAA$} satisfying $\AAAA(\xi_0)\oplus\GGG(\xi_0)=T_{\xi_0}\MMM$. Conversely, suppose that the Lie algebra of infinitesimal symmetries of $\Sigma=(f,g)$ possesses an abelian subalgebra $\AAAA=\vectR{v_1,\ldots,v_{n-m}}$. {By $\AAAA(\xi_0)\oplus\GGG(\xi_0)=T_{\xi_0}\MMM$ the vector fields $v_1, \ldots, v_{n-m}$ are linearly independent, so we} choose local coordinates $\tilde{\xi}=(\tilde{x},\tilde{w})$ such that $v_i=\vec{\tilde{x}_i}$, for $1\leq i\leq n-m$. In those coordinates, we have $g_j=A_j(\tilde{x},\tilde{w})\vec{\tilde{x}}+B_j(\tilde{x},\tilde{w})\vec{\tilde{w}}$. Since $\GGG$ is of constant rank $m$ and satisfies $\AAAA(\tilde{\xi}_0)\oplus\GGG(\tilde{\xi}_0)=T_{\tilde{\xi}_0}\MMM$, via a suitable feedback transformation {we choose generators of $\GGG$ as} $\tilde{g}_j=A_j(\tilde{x},\tilde{w})\vec{\tilde{x}}+\vec{\tilde{w}_j}${ (to simplify notations, we skip the "tildes" and denote $\tilde{g}_j$ by $g_j$)}. Using that $v_i$ are symmetries of $\Sigma$, that is $\lb{v_i}{g_j}\in\GGG$, we deduce that $A_j=A_j(\tilde{w})$, therefore we actually have $\lb{v_i}{g_j}=0$. Moreover, $\GGG$ is involutive so we deduce that $\lb{g_j}{g_k}=0$. {Therefore, all vector fields $v_i$ and $g_j$ commute and thus} there exist coordinates $\xi=(x,w)$ such that $v_i=\vec{x_i}$, for $1\leq i\leq n-m$, and $g_j=\vec{w_j}$, for $1\leq j\leq m$. The fields $v_i$ are symmetries of $\Sigma$ so $\lb{v_i}{f}\in\GGG$ implying that $f=F\vec{x}+\textsc{f}\vec{w}$, where $F=F(w)$ and we achieve $\textsc{f}=0$ by a suitable feedback transformation.
\end{enumerate}
\end{proof}
\subsection{Normal form of trivial systems possessing an almost abelian Lie subalgebra of infinitesimal symmetries}

\Cref{thm:charac_trivial_sys} of the previous subsection asserts that the Lie algebra of infinitesimal symmetries of trivialisable systems possesses an abelian subalgebra complementary to the distribution $\GGG$. In this subsection, {we study a particular case when the abelian subalgebra $\AAAA$ is actually an abelian ideal of codimension one of a subalgebra $\LLLL$ of symmetries, the latter acting transitively on $\MMM$.} We give a normal form of control-affine systems (with a, necessarily, scalar control) possessing such Lie algebra of symmetries.

\begin{definition}[Almost abelian Lie algebra] 
Let $\LLLL$ be a real Lie algebra{;} following the definition of \emph{\cite{burde2011Abelianidealsmaximal}}, we call $\LLLL$ \emph{almost abelian} if it has an abelian ideal $\IIII$ of codimension {one}.
\end{definition}
%
It is a simple application of Lie algebra homology to deduce that an almost abelian Lie algebra (possibly of infinite dimension) is isomorphic to the semi-direct product $\LLLL\cong\IIII\rtimes \vectR{v_0}$ and that its structure is determined by the action of $v_0$ on $\IIII$, namely by 
\begin{align*}
    \ad{v_0}{}\,:\,\IIII & \longrightarrow\IIII\\
    v&\longmapsto \lb{v_0}{v}.
\end{align*}
\noindent 
Moreover, two almost abelian Lie algebras $\LLLL=\IIII\rtimes\vectR{v_0}$ and $\tilde{\LLLL}=\tilde{\IIII}\rtimes\vectR{\tilde{v}_0}$ are isomorphic if and only if there exists a real invertible transformation $P:\IIII\rightarrow\tilde{\IIII}$ such that $P\ad{v_0}{}=\mu\,\ad{\tilde{v}_0}{}P$ for some $\mu\in\RR^*$; see \cite[Proposition 11]{avetisyan2016StructurealmostAbelian}. Therefore, isomorphism classes of almost abelian Lie algebras correspond to similarity classes of the linear operator $\ad{v_0}{}$ (up to multiplication by a scalar). In particular, if $\LLLL$ is finite dimensional, the similarity classes of $\ad{v_0}{}$ are given by the Jordan normal forms. In the following theorem, we consider the simplest case where $\ad{v_0}{}$ is diagonalisable over $\RR$ and give a normal form of control-affine systems that have $\LLLL$ as Lie subalgebra of infinitesimal symmetries.
\begin{theorem}[Almost abelian infinitesimal symmetries]\label{thm:normal_form_trivial_almost_abelian}
Consider a control-affine system $\Sigma=(f,g)${ on an $n$-dimensional state manifold $\MMM$ and with scalar control. Assume that the Lie algebra of infinitesimal symmetries possesses a Lie subalgebra $\LLLL$ for which the following conditions hold} at $\xi_0$: $f(\xi_0)\notin\GGG(\xi_0)$, $\LLLL$ acts transitively on $\MMM$, $\LLLL=\IIII\rtimes\vectR{v_0}$ is almost abelian, its abelian ideal satisfies $\IIII(\xi_0)\oplus\GGG(\xi_0)=T_{\xi_0}\MMM$, and the operator $\ad{v_0}{}$ is non-singular and diagonalisable over $\RR$. Then, $\Sigma$ is, locally around $\xi_0$, feedback equivalent to {a} {trivial system} {of the form}
\begin{align*}
    \Sigma_{\lambda}\,:\,\left\{\begin{array}{rl}
        \dot{x}_i &= \eta_i(w+1)^{\lambda_i}\quad 1\leq i\leq n-1 \\ 
        \dot{w} &= u
    \end{array}\right.,\quad u\in \RR
\end{align*}
\noindent
around $(x_0,0)\in\RR^{n-1}\times \RR$, where the $\lambda_i$'s are the eigenvalues of $\ad{v_0}{}$ and $\eta_i$ are constants equal $0$ or $1$, {with} $\eta_1=1$.
\end{theorem}
Observe that the assumptions: $\LLLL$ is almost abelian, acts transitively on $\MMM$, and satisfies $\IIII(\xi_0)\oplus\GGG(\xi_0)=T_{\xi_0}\MMM$, are quiet restrictive on the number of controls of $\Sigma$. Namely, they imply that $\rk{\GGG}=1$, i.e. $\Sigma$ has a scalar control. 
\begin{remark}[Converse implication]
    By a straightforward computation, we see that the vector fields 
    \begin{align*}
        v_i=\vec{x_i},\quad\forall\,1\leq i\leq n-1,\quad\textrm{and}\quad v_0=\sum_{i=1}^{n-1}\lambda_ix_i\vec{x_i}+(w+1)\vec{w} 
    \end{align*}
    \noindent
    are infinitesimal symmetries of $\Sigma_{\lambda}$. Thus the Lie algebra of infinitesimal symmetries of $\Sigma_{\lambda}$ possesses an almost abelian subalgebra $\LLLL=\vectR{v_1,\ldots,v_{n-1},v_0}$. Therefore, the theorem is actually an "if and only if" statement since all other assumptions are feedback invariant. 
\end{remark}
\begin{proof}
    Consider the control-affine system $\Sigma=(f,g)$ given by vector fields $f$ and $g$, and let $n$ vector fields $v_1,\ldots,v_{n-1},v_0$ generate the $n$-dimensional Lie subalgebra $\LLLL=\vectR{v_1,\ldots,v_{n-1},v_0}$ of {the algebra of} infinitesimal symmetries, which by assumption is almost abelian, whose abelian ideal is $\IIII=\vectR{v_1,\ldots,v_{n-1}}$ and since $\ad{v_0}{}$ is diagonalisable over $\RR$ we {conclude} that
    \begin{align}\label{eq:mult_table_llll}
        \lb{v_0}{v_i}=\lambda_iv_i,\quad\forall\,1\leq i\leq n-1.
    \end{align}
    \noindent
    By statement \cref{thm:charac_trivial_sys:2} of \cref{thm:charac_trivial_sys}, $\Sigma$ is locally trivialisable and following the proof of that theorem we deduce that there exists local coordinates $(x,\hat{w})$ around $(x_0,0)$ such that $v_i=\vec{x_i}$, for $1\leq i\leq n-1$, $g=\vec{\hat{w}}$, and $f=\sum_{i=1}^{n-1}F_i\vec{x_i}$, where $F_i=F_i(\hat{w})$. {Express the infinitesimal symmetry} $v_0=\sum_{i=1}^{n-1}\gamma_i\vec{x_i} +\bar{\delta}\vec{\hat{w}}$, where $\gamma_i=\gamma_i(x)$ since $v_0$ is a symmetry of $\GGG=\distrib{\vec{\hat{w}}}$, and $\bar{\delta}=\bar{\delta}(x,\hat{w})$. Using the commutativity relations \cref{eq:mult_table_llll} that have not been changed by applying diffeomorphisms, we obtain
    \begin{align*}
        v_0 &= -\sum_{i=1}^{n-1}\left(\lambda_ix_i+a_i\right)\vec{x_i}+\bar{\delta}(\hat{w})\vec{\hat{w}},\quad a_i\in\RR.
    \end{align*}
    \noindent
    To avoid unnecessary computations, we replace $v_0$ by$-v_0-\sum_{i=1}^{n-1}a_iv_i\in\LLLL${, having the same properties}, so we can assume $v_0=\sum_{i=1}^{n-1}\lambda_ix_i\vec{x_i}+\delta(\hat{w})\vec{\hat{w}}$, where $\delta(\hat{w})=-\bar{\delta}(\hat{w})$. Using the fact that $v_0$ is a symmetry of $f$, i.e. $\lb{v_0}{f}\in\GGG$, we deduce the following equations: 
    \begin{align}\label{eq:sys_sym_f}
        \delta(\hat{w})\frac{\diff F_i}{\diff \hat{w}}(\hat{w})-\lambda_iF_i(\hat{w})=0,\quad \forall\,1\leq i\leq n-1.
    \end{align}
    \noindent
    {The assumption} $f(\xi_0)\notin\GGG(\xi_0)${implies} that $(F_1,\ldots,F_{n-1})(0)\neq0\in\RR^{n-1}$, {so} suppose that $F_1(0)\neq 0$ (renumber the $x_i$'s if necessary) and thus $F_1(\hat{w})=c+h(\hat{w})$ with $c=F_1(0)$ for {a function $h$ satisfying} $h(0)=0$. Replacing $x_1$ by $\frac{x_1}{c}$ we may assume that $F_1(\hat{w})=1+h(\hat{w})$. Moreover, by \cref{eq:sys_sym_f} we conclude that $h'(0)\neq0$; recall that by assumption $\ad{v_0}{}$ is {non-}singular {and} hence $\lambda_i\neq0$ for all $1\leq i\leq n-1$. Thus, ${w}=(1+h(\hat{w}))^{1/\lambda_1}-1$ is a diffeomorphism around $\hat{w}_0=0$ and, to simplify notations, we keep the symbols ${f}$ and ${v}_0$ for {those} vector fields {represented using the coordinate $w$}. {We have $F_1(w)=(1+w)^{\lambda_1}$, thus relation} \cref{eq:sys_sym_f} gives{, with $\hat{w}$ renamed $w$,} ${\delta}({w})={w}+1$ and implies that the functions ${F}_2,\ldots,{F}_{n-1}$ satisfy 
    \begin{align*}
        (1+{w})\frac{\diff {F}_i}{\diff {w}}(w)=\lambda_i{F}_i({w}),\quad \textrm{for}\quad 2\leq i\leq n-1.
    \end{align*}
    \noindent
    Solving those equations around ${w}_0=0$ gives ${F}_i=c_i\left({w}+1\right)^{\lambda_i}$, with $c_i\in\RR$. 
    Thus, normalising $x_i$ with the non-zero $c_i$ we get that $\Sigma$ takes the form of $\Sigma_{\lambda}$ around $(x_0,0)\in\RR^n$.
\end{proof}
The above theorem generalises our previous results on the Lie algebra of infinitesimal symmetries of {hyperbolic and }parabolic systems for which we have {$n=3$ and, respectively, $(\lambda_1,\lambda_2)=(1,-1)$ and $(\lambda_1,\lambda_2)=(2,1)$, see  \cite{SCHMODERER2022105397}.} 

The situation is much more involved when $f(\xi_0)\in\GGG(\xi_0)$ as the following proposition illustrates. 
\begin{proposition}\label{prop:normal_form_trivial_system_non_generic}
Consider a control-affine system $\Sigma=(f,g)$ with a $n$-dimensional state manifold and a scalar control, and let $\LLLL$ be a Lie subalgebra of infinitesimal symmetries. Suppose that the following conditions hold at $\xi_0$: $f(\xi_0)\in\GGG(\xi_0)$, there exists $k\geq1$ the smallest integer such that $g\wedge\adk{g}{k}{f}(\xi_0)\neq0$, $\LLLL$ acts transitively on $\MMM$, {$\LLLL\cong\IIII\rtimes \vectR{v_0}$} is almost abelian, $\IIII(\xi_0)\oplus\GGG(\xi_0)=T_{\xi_0}\MMM$, and the action {of} $\ad{v_0}{}$ is non-singular and diagonalisable over $\RR$. {Then, } the eigenvalues $\lambda_i$ of $\ad{v_0}{}$ {are such that} $\frac{\lambda_i}{\lambda_1}k$ are {positive }integers greater or equal than $k$, where $\lambda_1$ is the smallest{, in absolute value,} eigenvalue of $\ad{v_0}{}$. {Moreover, the system} $\Sigma$ is locally feedback equivalent to 
\begin{align*}
    \Sigma_{\lambda}^{0,k}\,:\,\left\{\begin{array}{rl}
        \dot{x}_1 &= w^k \\
        \dot{x}_i &= \eta_i(w^k)^{\lambda_i/\lambda_1}\quad 2\leq i\leq n-1 \\ 
        \dot{w} &= u
    \end{array}\right.,\quad u\in \RR,
\end{align*}
\noindent
around $(x_0,0)\in\RR^n$, where $\eta_i$ are constants equal to $0$ or $1$. 
\end{proposition}
Observe that the normal form $\Sigma_{\lambda}^{0,k}$ defines {a} polynomial system {since $k\frac{\lambda_i}{\lambda_1}$ are {positive} integers}. {Moreover, although the eigenvalues $\lambda_i$ can be any, the conditions of the above proposition imply severe restrictions on them. Namely, all $\lambda_i$'s have the same sign and $\lambda_i=\frac{l_i}{k}\lambda_1$, where the $l_i$'s are integers satisfying $l_i\geq k$.} Furthermore, it is a classical fact that (under the above assumptions) the integer $k$ is an invariant of feedback transformations. Hence if $k\neq k'$, then $\Sigma_{\lambda}^{0,k}$ and $\Sigma_{\lambda}^{0,k'}$ are not locally feedback invariant around $(x_0,0)$.
\begin{proof}
The beginning of the proof is the same as that of the previous theorem up to equation \cref{eq:sys_sym_f}, so we start from there. {By} $f(\xi_0)\in\GGG(\xi_0)$ {we have} $(F_1,\ldots,F_{n-1})(0)=0\in\RR^{n-1}$ and {due to} $g\wedge\adk{g}{k}{f}(0)\neq0$ we assume that $\frac{\diff^k F_1}{\diff \hat{w}^k}(0)\neq0$, if necessary relabel the $x_i$-coordinates. {Using the Taylor expansion we can write $F_1=\hat{w}^kH(\hat{w})$, where $H(0)\neq0$ and we can suppose $H(0)>0$, if not, replace $x_1$ by $-x_1$. We apply the local diffeomorphism ${w}=\hat{w}\left(H(\hat{w})\right)^{1/k}$}, that maps $f$ and $v_0$ into vector fields that, for simplicity, we denote again by ${f}$ and ${v}_0$, respectively. We have ${F}_1({w})={w}^k$, so equation \cref{eq:sys_sym_f}, expressed in the $w$-coordinate, implies, for $i=1$, that ${\delta}({w})=\frac{\lambda_1}{k}{w}$ and
\begin{align*}
    \lambda_1{w}\frac{\diff {F}_i}{\diff {w}}(w)=\lambda_i k {F}_i({w}),\quad\textrm{for}\quad 2\leq i\leq n-1.
\end{align*}
\noindent
Solving those equations around ${w}_0=0$ {implies $|F_i|=c_i|w_i|^{k\lambda_i/\lambda_1}$}, with $c_i\in\RR$. The only $C^{\infty}$-solutions are those given by either $c_i=0$ or by $c_i\neq0$ with $\frac{\lambda_i}{\lambda_1}k$ being a {positive} integer{. The corresponding smooth solution is $F_i=c_iw^{k\lambda_i/\lambda_1}$ and by definition of $k$ we have $k\frac{\lambda_i}{\lambda_1}\geq k$ {(otherwise, $g\wedge\adk{f}{k'}{g}(\xi_0)\neq0$, where $k'=\frac{\lambda_i}{\lambda_1}k<k$, contradicting the definition of $k$)}. Thus $|\lambda_i|\geq|\lambda_1|$ and $\lambda_1$ is, indeed, the eigenvalue of $\ad{v_0}{}$ of minimal absolute value.}. 
Finally, normalising the coordinates $x_i$, with $c_i\neq0$, we obtain the desired form $\Sigma_{\lambda}^{0,k}$.
\end{proof}
The previous proposition describes all smooth systems having an almost abelian Lie subalgebra of symmetries for which $k$ exists, in particular all analytic systems. Notice that for a single-input analytic system either $k$ exists or, if not, then it is locally feedback equivalent to a trivial system $\dot{x}=c,\,\dot{w}=u$, where $c\in\RR^{n-1}$. In the $C^{\infty}$ category there are, however, systems for which $k$ does not exist but the symmetry algebra possesses an almost abelian subalgebra. For example, consider around $(x_0,0)$ the system
\begin{align*}
    \left\{\begin{array}{rl}
        \dot{x}_1 &= \mathfrak{f}(w) \\
        \dot{x}_i &= \mathfrak{f}(w)^{\lambda_i/\lambda_1}\quad 2\leq i\leq n-1 \\ 
        \dot{w} &= u
    \end{array}\right.,\quad\textrm{with}\quad \mathfrak{f}(w)=\exp\left(-\frac{1}{w^2}\right),\;\mathfrak{f}(0)=0,
\end{align*}
\noindent
and $\lambda_i/\lambda_1\in\NN^*$. By a straightforward calculation, one may check that the system possesses an almost abelian Lie subalgebra of infinitesimal symmetries but, obviously, $k$ does not exist at $(x_0,0)$ and thus it is not feedback equivalent to $\Sigma_\lambda^{0,k}$.
\section{Trivial systems on 3D-manifolds}\label{sec:trivial_systems_3d}
In this section, we study trivial system on $3$-dimensional manifolds with scalar control. Our aim is to give a new version of the results of \cite{serres2009Controlsystemszero} and to extend them by giving several normal forms. \\

\noindent
Throughout this section, we consider a control-affine system $\Sigma=(f,g)$ of the form
\begin{align*}
    \Sigma\,:\,\dot{\xi}=f(\xi)+g(\xi)u,\quad u\in\RR,
\end{align*}
\noindent
where the state $\xi$ belongs to a $3$-dimensional manifold {$\MMM$} and the vector fields $f$ and $g$ satisfy{, at any $\xi\in\MMM$,} the following regularity assumptions
\begin{enumerate}[label=\textrm{(A\arabic*)},ref=\textrm{(A\arabic*)},font=\normalfont]
    \item\label{assumption:1} $f\wedge g\wedge\lb{g}{f}\neq0$, 
    \item\label{assumption:2} $g\wedge\lb{g}{f}\wedge\lb{g}{\lb{g}{f}}\neq0$.
\end{enumerate}
To any control-affine system $\Sigma=(f,g)$ we attach $6$ structure functions uniquely defined by the following decompositions: 
\begin{align}\label{eq:structure_functions_definition}
    \begin{array}{rl}
        \lb{f}{\lb{f}{g}} &= k_1g+k_2\lb{g}{f}+k_3\lb{g}{\lb{g}{f}},\\
        \lb{g}{\lb{g}{f}} &=\lambda_1 f+\lambda_2g+\lambda_3\lb{g}{f}.
    \end{array}
\end{align}
\noindent
Observe that assumption \cref{assumption:2} implies that $\lambda_1\neq0$. We now define two different classes of pairs $(f,g)$.
\begin{definition}[(Semi)-canonical pairs]\label{def:pair_control_affine_system}
We call the pair $(f,g)$ \emph{semi-canonical} if $k_3\equiv 0$, and we will denote it by $(f_s,g)$. If, additionally $\lambda_1\equiv \pm1$, then we call $(f,g)$ a \emph{canonical} pair and we denote it by $(f_c,g_c)$.
\end{definition}
Observe that a semi-canonical pair is characterised by the inclusion $\lb{f}{\lb{f}{g}}\in\distrib{g,\lb{g}{f}}$, which is the property of the singular vector field of $\Sigma$ (justifying the notation $(f_s,g)$ for a semi-canonical pair). Under assumption \cref{assumption:2}, that singular vector field is unique and can be computed using the singular control, thus we are not surprised that the following proposition shows that a semi-canonical pair exists (but observe that our proof does not require {using} the machinery of singular controls). 

The following proposition shows that a control-affine system is always feedback equivalent to a system given by a semi-canonical and even by a canonical pair. Moreover, those pairs can be \emph{explicitly} constructed, meaning that the feedback transformations bringing $(f,g)$ into $(f_s,g)$ or into $(f_c,g_c)$ are constructed with {differentiation and} algebraic operations only (no differential equations to be solved). Furthermore, a canonical pair is \emph{unique} (up to {multiplying $g_c$ by $-1$}) hence its structure functions are feedback {equivariants} (up to a discrete involution, see \cref{rk:unique_canonical_pair} below). 
\begin{proposition}[Existence of semi-canonical and canonical pairs]\label{prop:existence_pairs}
Consider a control-affine system $\Sigma=(f,g)$ {satisfying assumptions \cref{assumption:1} and \cref{assumption:2}}. Then, the following statements hold globally: 
\begin{enumerate}[label=(\roman*),ref=\textit{(\roman*)}]
    \item \label{prop:existence_pairs:1} $\Sigma$ is globally feedback equivalent to $\Sigma_s=(f_s,g)$, where $(f_s,g)$ is a semi-canonical pair;
    \item \label{prop:existence_pairs:2} $\Sigma$ is globally feedback equivalent to $\Sigma_c=(f_c,g_c)$, where $(f_c,g_c)$ is a canonical pair.
\end{enumerate}
\noindent
Moreover $(f_s,g)$ and $(f_c,g_c)$ can be explicitly constructed via the following feedback transformations
\begin{align*}
    f_s=f+gk_3,\quad f_c=f_s,\quad \textrm{and}\quad g_c=\left|\lambda_1\right|^{-1/2}g.
\end{align*}
\noindent
Furthermore, the canonical pair $(f_c,g_c)$ is unique up to $g_c\mapsto -g_c$. 
\end{proposition}
{The proof of the above proposition is based on the following lemma, which gives some relations between the structure functions {$k_1$, $k_2$, $k_3$, $\lambda_1$, $\lambda_2$, and $\lambda_3$,} and shows how they {change} under feedback transformations. Its proof is a straightforward computation that we detail in \cref{apdx:transform_structure_functions}.}
\begin{lemma}\label{lem:transform_structure_functions}
Consider a control-affine system $\Sigma=(f,g)$ with structure functions $(k_1,k_2,k_3)$ and $(\lambda_1,\lambda_2,\lambda_3)$. Then, the following relations hold: 
\begin{subequations}
\begin{align}
  	\dL{f}{\lambda_1}  &= - k_2\lambda_1-\dL{g}{\lambda_1 k_3}, \label{eq:relation_struct_functions:1} \\ 
    \dL{f}{\lambda_2}-\lambda_3 k_1 +\dL{g}{k_1}&= - k_2\lambda_2-\dL{g}{\lambda_2 k_3}, \label{eq:relation_struct_functions:2} \\ 
    \dL{f}{\lambda_3}-\lambda_2 &=-k_3\lambda_1-\dL{g}{k_2}-\dL{g}{\lambda_3 k_3}. \label{eq:relation_struct_functions:3}
\end{align}
\end{subequations}
Under a feedback transformation $f\mapsto\tilde{f}=f+g\alpha$ and $g\mapsto \tilde{g}=g\beta$, where $\alpha$ and $\beta$ are smooth scalar functions satisfying $\beta\neq0$, the structure functions are transformed by 
\begin{align}
    &\begin{array}{rl}
        \tilde{k}_1&= k_1+\dL{\lb{g}{f}}{\alpha} + \frac{1}{\beta}\left(\dL{f}{\gamma}+\alpha\dL{g}{\gamma}-\gamma\dL{g}{\alpha}\right)+\tilde{k}_2\frac{\gamma}{\beta}\\ 
        &\quad+\tilde{k}_3\left(\dL{\lb{g}{f}}{\beta}+\dL{g}{\gamma}-\gamma\dL{g}{\ln(\beta)}\right),\\
        \tilde{k}_2&=k_2-\dL{f}{\ln(\beta)}-\frac{\gamma}{\beta}-\alpha\dL{g}{\ln(\beta)}+\tilde{k}_3\dL{g}{\beta},\\
        \tilde{k}_3&=\frac{1}{\beta}\left( k_3-\alpha\right),
    \end{array} \label{eq:transform_struct_funct_k} \\
    \tilde{\lambda}_1&= \beta^2\lambda_1, \quad\tilde{\lambda}_2=\beta\lambda_2-\beta\lambda_1\alpha+\gamma\lambda_3  -\dL{\lb{g}{f}}{\beta}-\dL{g}{\gamma}+2\gamma\dL{g}{\ln(\beta)}, \quad\tilde{\lambda}_3=\beta\lambda_3+\dL{g}{\beta},\label{eq:transform_struct_funct_lambda}
\end{align}
\noindent
where $\gamma=\dL{f}{\beta}+\alpha\dL{g}{\beta}-\beta\dL{g}{\alpha}$.
\end{lemma}
\begin{proof}[Proof of \cref{prop:existence_pairs}.]
Consider $\Sigma=(f,g)$ {whose} structure functions {are} $(k_1,k_2,k_3)$ and $(\lambda_1,\lambda_2,\lambda_3)$. By \cref{eq:transform_struct_funct_k,eq:transform_struct_funct_lambda} we have that under feedback transformations $\beta\tilde{k}_3=k_3-\alpha$ and $\tilde{\lambda}_1=\beta^2\lambda_1$. Hence, choosing  $\alpha=k_3$ we obtain that the transformed pair $(\tilde{f},g)${, where $\tilde{f}=f+gk_3$,} is semi-canonical. Moreover, additionally choosing $\beta=\left|\lambda_1\right|^{-1/2}$, recall that $\lambda_1\neq0$ from assumption \cref{assumption:2}, yields a canonical pair $(f_c,g_c)=(\tilde{f},\tilde{g})${, where $\tilde{g}=\beta g$}. Clearly, the singular vector field $f_s=f_c$ is uniquely defined, but the canonical vector field $g_c$ is unique up to $g_c\mapsto \pm g_c$.
\end{proof}
Observe that for the canonical pair $(f_c,g_c)$ we additionally have $k_2\equiv 0$, due to \cref{eq:relation_struct_functions:1}. Thus, the canonical pair $(f_c,g_c)$ satisfies the following decomposition (renaming $k_1$ to $\kappa$, $\lambda_1$ to $\varepsilon$, $\lambda_2$ to $\mu$, $\lambda_3$ to $\nu$)
\begin{align}\tag{\ref*{eq:structure_functions_definition}'}\label{eq:canonical_structure_constants}
    \begin{array}{rl}
        \lb{f_c}{\lb{f_c}{g_c}} &= \kappa g_c,\\ 
        \lb{g_c}{\lb{g_c}{f_c}} &= \varepsilon f_c +\mu g_c+\nu \lb{g_c}{f_c},
    \end{array}
\end{align}
\noindent
where $\varepsilon=\pm1$. {Moreover, using \cref{eq:relation_struct_functions:2,eq:relation_struct_functions:3} we deduce that $\kappa$, $\mu$, and $\nu$ are related by
\begin{align}
    \dL{f_c}{\mu}-\nu\kappa+\dL{g_c}{\kappa}&=0 \tag{\ref*{eq:relation_struct_functions:2}'}\label{eq:relation_struct_functions:2_p}\\ 
    \dL{f_c}{\nu}-\mu&=0,\tag{\ref*{eq:relation_struct_functions:3}'}\label{eq:relation_struct_functions:3_p}
\end{align}
\noindent
from which we deduce {that the feedback invariants $\kappa$ and $\nu$ are associated via}
\begin{align}\label{eq:relation_kappa_nu}
    \dLk{f_c}{2}{\nu}-\nu\kappa+\dL{g_c}{\kappa}&=0.
\end{align}
\noindent
}Therefore, a canonical pair identifies explicitly a discrete invariant $\varepsilon=\pm1$ and  two constructible feedback invariant functions $\kappa$ and $\nu$ called, respectively, the \emph{curvature} and the \emph{centro-affine curvature} by analogy with Serres' work \cite{serres2007curvaturefeedbackclassification}{; see also \cite{wilkens1998Centroaffinegeometryplane}}. Observe that due to \cref{eq:relation_struct_functions:3_p} {above}, $\mu$ is determined by $\nu$. {Moreover, the curvature $\kappa$ determines the centro-affine curvature $\nu$ up to a affine part; i.e. if $f_c$ is rectified on $\vec{x}$, then $\nu$ is determined by $\kappa$ via \cref{eq:relation_kappa_nu} up to two functions $\nu_1$ and $\nu_0$ satisfying $\dL{f_c}{\nu_i}=0$.}
\begin{remark}\label{rk:unique_canonical_pair}
    A canonical pair is unique up to $g_c\mapsto -g_c$. Hence the centro-affine curvature $\nu$ is a feedback {equivariant} up to the involution $\nu\mapsto-\nu$ (which does not influence our conditions below). We will {get} back {to that} subtlety in \cref{prop:normal_forms}{,} {where} we will construct several normal forms. {On} the other hand, the curvature $\kappa$ is a true {feedback} invariant {(actually, a feedback equivariant that changes as $\phi^*\kappa$ under a diffeomorphism $\phi$).}
\end{remark}
\noindent
For a given control-affine system $\Sigma_c=(f_c,g_c)$, given by a canonical pair, we will denote by { $(\varepsilon,\kappa,\nu)$ the triple of invariants}. Although the canonical pair can be constructed without much work, for the sake of completeness, we give the expression of $(\varepsilon,\kappa,\nu)$ for an arbitrary control-affine system. In particular, observe that our formula for the curvature $\kappa$ generalises the one given in \cite[p. 376]{agrachev2013ControlTheoryGeometric}, where $k_3$ is already normalised to $0$.

\begin{proposition}[Invariants of control-affine systems]\label{prop:invariants_control_systems}
    Consider a control-affine system $\Sigma=(f,g)$ {on}  a $3$-dimensional {state-space} manifold, and with scalar control, and let $(k_1,k_2,k_3)$ and $(\lambda_1,\lambda_2,\lambda_3)$ be structure functions defined by \cref{eq:structure_functions_definition}. Then, the invariants $(\varepsilon,\kappa,\nu)$ are given by: 
    \begin{align}
        \varepsilon=\sgn{\lambda_1},\quad\kappa&=k_1+\frac{1}{2}\dL{f}{k_2-\dL{g}{k_3}}+\frac{1}{4}\left(k_2-\dL{g}{k_3}\right)^2+\dL{\lb{g}{f}}{k_3}+\frac{1}{2}k_3\dL{g}{k_2-\dL{g}{k_3}},\nonumber\\
        \textrm{and}\quad\nu&=\left|\lambda_1\right|^{-1/2}\left(\lambda_3-\frac{1}{2}\dL{g}{\ln\left|\lambda_1\right|}\right). \label{eq:formula_curavtures_struct_functions}
    \end{align}
    \noindent
\end{proposition}
Our formula for the curvature $\kappa$ is, indeed, a generalisation of {that in} \cite{agrachev2013ControlTheoryGeometric} because if $k_3=0$ (i.e. we suppose that $f$ is the singular vector field {$f_s$}), then \cref{eq:formula_curavtures_struct_functions} reads 
\begin{align*}
    \kappa &= k_1+\frac{1}{2}\dL{f}{k_2}+\frac{1}{4}\left(k_2\right)^2,
\end{align*}
\noindent
that is to say, exactly as the formula given by Agrachev \cite{agrachev1998FeedbackInvariantOptimalControl,agrachev2013ControlTheoryGeometric}. 
\begin{proof}
    Consider a pair $(f,g)$ with structure functions $(k_1,k_2,k_3)$ and $(\lambda_1,\lambda_2,\lambda_3)$. 
    To deduce the {expression} of the invariants $(\varepsilon,\kappa,\nu)$, we apply \cref{eq:transform_struct_funct_k,eq:transform_struct_funct_lambda} with $\alpha=k_3$ and $\beta=\left|\lambda_1\right|^{-1/2}$, namely the feedback transformation that constructs the canonical pair. We detail the computation for $\kappa$ {and left the computation for $\nu$ to the reader}. Recall, that applying a feedback to construct the canonical pair, we obtain as a by-product $\tilde{k}_2=0$ {(see the proof of \cref{prop:existence_pairs})}. First we have,
    \begin{align*}
        \gamma = -\frac{1}{2}\left|\lambda_1\right|^{-1/2}\dL{f}{\Lambda}-\frac{1}{2}k_3\left|\lambda_1\right|^{-1/2}\dL{g}{\Lambda}-\left|\lambda_1\right|^{-1/2}\dL{g}{k_3} = \left|\lambda_1\right|^{-1/2}\left(-\frac{1}{2}\dL{f}{\Lambda}-\frac{1}{2}k_3\dL{g}{\Lambda}-\dL{g}{k_3}\right),
    \end{align*}
    \noindent
    where $\Lambda=\ln\left|\lambda_1\right|$. Second, using $\tilde{k}_2=0$ we deduce
    \begin{align*}
        k_2&=-\frac{1}{2}\dL{f}{\Lambda}-\frac{1}{2} k_3\dL{g}{\Lambda}+\left|\lambda\right|^{1/2}\gamma =-\dL{f}{\Lambda}-k_3\dL{g}{\Lambda}-\dL{g}{k_3}.
    \end{align*}
    \noindent
    Therefore, inserting the last expression of $k_2$ into $\gamma$, we get $\gamma=\frac{1}{2}\left|\lambda_1\right|^{-1/2}\left(k_2-\dL{g}{k_3}\right)$. Now the curvature reads $\kappa=\tilde{k}_1$, i.e 
    \begin{align*}
        \kappa &= k_1+\dL{\lb{g}{f}}{k_3}+\left|\lambda_1\right|^{1/2}\left(\frac{1}{2}\dL{f}{\left|\lambda_1\right|^{-1/2}\left(k_2-\dL{g}{k_3}\right)}+\frac{1}{2}k_3\dL{g}{\left|\lambda_1\right|^{-1/2}\left(k_2-\dL{g}{k_3}\right)}\right.\\
        &\quad \left.-\frac{1}{2}\left|\lambda_1\right|^{-1/2}\left(k_2-\dL{g}{k_3}\right)\dL{g}{k_3}\right),\\ 
        &=k_1+\dL{\lb{g}{f}}{k_3}-\frac{1}{2}\left(k_2-\dL{g}{k_3}\right)\dL{g}{k_3}\\
        &\quad+\frac{1}{2}\left|\lambda_1\right|^{1/2}\left(\left|\lambda_1\right|^{-1/2}\dL{f}{k_2-\dL{g}{k_3}}-\frac{1}{2}(k_2-\dL{g}{k_3})\left|\lambda_1\right|^{-3/2}\dL{f}{\left|\lambda_1\right|}\right) \\
        &\quad+ \frac{1}{2}\left|\lambda_1\right|^{1/2}k_3\left(\left|\lambda_1\right|^{-1/2}\dL{g}{k_2-\dL{g}{k_3}}-\frac{1}{2}(k_2-\dL{g}{k_3})\left|\lambda_1\right|^{-3/2}\dL{g}{\left|\lambda_1\right|}\right)\\
        &=k_1+\dL{\lb{g}{f}}{k_3}-\frac{1}{2}\left(k_2-\dL{g}{k_3}\right)\dL{g}{k_3} + \frac{1}{2}\dL{f}{k_2-\dL{g}{k_3}}\\
        &\quad-\frac{1}{4}\left(k_2-\dL{g}{k_3}\right)\dL{f}{\Lambda} + \frac{1}{2}k_3\dL{g}{k_2-\dL{g}{k_3}}-\frac{1}{4}\left(k_2-\dL{g}{k_3}\right)k_3\dL{g}{\Lambda} \\ 
        &= k_1+\dL{\lb{g}{f}}{k_3} + \frac{1}{2}\dL{f}{k_2-\dL{g}{k_3}}+\frac{1}{2}k_3\dL{g}{k_2-\dL{g}{k_3}}\\
        &\quad-\frac{1}{4}\left(k_2-\dL{g}{k_3}\right)\left(\dL{f}{\Lambda}+\dL{g}{\Lambda}+2\dL{g}{k_3}\right) \\
        &= k_1+\frac{1}{2}\dL{f}{k_2-\dL{g}{k_3}}+\frac{1}{4}\left(k_2-\dL{g}{k_3}\right)^2+\dL{\lb{g}{f}}{k_3}+\frac{1}{2}k_3\dL{g}{k_2-\dL{g}{k_3}}.
    \end{align*}
\end{proof}
\noindent
Now consider a trivial system{,} whose state $(x,y,w)$ belongs to a $3$-dimensional manifold {$\MMM$,}
\begin{align*}
    (T)\,:\,\left\{\begin{array}{rl}
        \dot{x} &= F_1(w) \\
        \dot{y} &= F_2(w) \\
        \dot{w} &= u
    \end{array}\right.,\quad (x,y,w)\in\MMM,\quad u\in\RR.
\end{align*}
\noindent
Notice that $(T)$ is, in general, not given by a canonical pair but is given by a semi-canonical pair since $\lb{f}{\lb{f}{g}}=0$. Clearly, for trivial systems we have $\kappa=0$, but the converse is not true as discovered in \cite{serres2009Controlsystemszero} and as we will show in the following theorem. 
\begin{theorem}[Characterisation of trivial systems]\label{thm:characterisation_trivialisable}
Consider a control-affine system {$\Sigma=(f,g)$ together with its structure functions $\kappa$ and $\nu$. Then, $\Sigma$ is locally trivialisable if and only if its canonical form} $\Sigma_c=(f_c,g_c)$ {satisfies} 
\begin{align}\label{eq:conditions_trivialisable}
    \kappa=0,\quad\dL{f_c}{\nu}=0,\quad\textrm{and}\quad\dL{\lb{f_c}{g_c}}{\nu}=0.
\end{align}
\end{theorem}
\noindent
Observe that the conditions of \cref{eq:conditions_trivialisable} can explicitly be tested on {the} control-affine system $\Sigma=(f,g)$. Indeed, with the help of \cref{prop:existence_pairs}, we explicitly construct the canonical pair $(f_c,g_c)$ of $\Sigma$ for which the invariants $\kappa$ and $\nu$ can be computed by algebraic operations only. Another way to test condition \cref{eq:conditions_trivialisable} on an arbitrary control-affine system {$\Sigma=(f,g)$} is to compute the invariants $\kappa$ and $\nu$ using \cref{eq:formula_curavtures_struct_functions} and then to evaluate \cref{eq:conditions_trivialisable} with $f_c=f+gk_3$ and $g_c=\left|\lambda_1\right|^{-1/2}g$. 
\begin{proof}
We begin with necessity {and} suppose that $\Sigma$ is trivialisable. Then, $(T)$ is given by $f=F_1(w)\vec{x}+F_2(w)\vec{y}$ and $g=\vec{w}$ (which, a priori, is not a canonical pair), whose structure functions are $k_1=k_2=k_3=0$, $\lambda_1=\lambda_1(w)$, $\lambda_2=0$, and $\lambda_3=\lambda_3(w)$. In particular, observe that $\lambda_1$ and $\lambda_3$ satisfy $\dL{f}{\lambda_1}=\dL{\lb{g}{f}}{\lambda_1}=0$ and $\dL{f}{\lambda_3}=\dL{\lb{g}{f}}{\lambda_3}=0$. {As in the proof of \cref{prop:existence_pairs},} to transform the pair $(f,g)$ of $(T)$ into the canonical pair {$(f_c,g_c)$} we use $\beta=\left|\lambda_1\right|^{-1/2}$, which therefore satisfies $\dL{f}{\beta}=0$ and $\dL{\lb{g}{f}}{\beta}=0$. {Now, using \cref{eq:transform_struct_funct_lambda,eq:transform_struct_funct_k} of \cref{lem:transform_structure_functions}, we calculate the structure functions of $(f_c,g_c)=(f,g\beta)$ {which are} $\tilde{\kappa}=0$, $\tilde{\varepsilon}=\pm1$, $\tilde{\mu}=0$, and $\tilde{\nu}=\tilde{\lambda}_3=\beta\lambda_3+\dL{g}{\beta}$.} 
Hence, for the canonical pair {$(f_c,g_c)$} of $(T)$ we have 
\begin{align*}
    \dL{f_c}{\tilde{\nu}} &= \beta\dL{f}{\lambda_3}+\dL{f}{\dL{g}{\beta}}= \beta\dL{f}{\lambda_3}+\dL{g}{\dL{f}{\beta}}-\dL{\lb{g}{f}}{\beta}= 0, \\ 
    \dL{\lb{g_c}{f_c}}{\tilde{\nu}} &= \beta\dL{\lb{g}{f}}{\tilde{\nu}}=\beta\left(\dL{\lb{g}{f}}{\lambda_3}+\dL{\lb{g}{f}}{\dL{g}{\beta}}\right) = \beta\left(\dL{g}{\dL{\lb{g}{f}}{\beta}}-\dL{\lb{g}{\lb{g}{f}}}{\beta} \right)\\ 
    & = \beta\left(-\lambda_1\dL{f}{\beta}-\lambda_3\dL{\lb{g}{f}}{\beta} \right)=0,
\end{align*}
\noindent
and the necessity of \cref{eq:conditions_trivialisable} is proved.

Now, conversely, suppose that $\Sigma_c$, given by its canonical pair $(f_c,g_c)$, satisfies \cref{eq:conditions_trivialisable}. {First, {due} to \cref{apdx:exists-diffeo:lemma} of \cref{apdx:exists-diffeo}, we apply a diffeomorphism $({x},{y},{w})=\phi(\xi)$ {that simultaneously rectifies the distribution $\FFF=\distrib{f_c,\lb{g_c}{f_c}}$ and the vector field $g_c$, that is} $\phi_*\FFF=\distrib{\vec{{x}},\vec{{y}}}$ and $\phi_* g_c=\vec{{w}}$. In those coordinates, we have $f_c=f_1\vec{{x}}+f_2\vec{{y}}$, with $f_i=f_i({x},{y},{w})$, and we have $\nu=\nu({w})$ since $\dL{f_c}{\nu}=\dL{\lb{g_c}{f_c}}{\nu}=0$ and $f_c\wedge g_c\wedge\lb{g_c}{f_c}\neq0$. Therefore, using relation \cref{eq:canonical_structure_constants} we deduce that $f_c$ satisfies the following two equations (notice that equation \cref{eq:relation_struct_functions:3_p} together with $\dL{f_c}{\nu}=0$ imply that $\mu=0$)
\begin{align*}
    \lb{f_c}{\lb{f_c}{g_c}} &= 0\quad\textrm{and}\quad\lb{g_c}{\lb{g_c}{f_c}}=\varepsilon f_c + \nu({w})\lb{g_c}{f_c}.
\end{align*}
\noindent
The second equation reads
\begin{align}\label{eq:ode_trivial_f}
    \frac{\partial^2f_c}{\partial{w}^2}=\varepsilon f_c+\nu({w})\frac{\partial f_c}{\partial{w}},
\end{align}
\noindent
and, interpreted as a second order linear ODE with respect to ${w}$ and with parameters $({x},{y})$, admits local solutions of the form 
\begin{align}\label{eq:ode_trivial_f_solution}
    f_c({x},{y},{w})=F_1({w})\begin{pmatrix}
{a}_1 \\ {a}_2 \\0
\end{pmatrix} + F_2({w})\begin{pmatrix}
{b}_1 \\ {b}_2 \\0
\end{pmatrix} = F_1({w}) {A} + F_2({w}){B}.
\end{align}
\noindent
In \cref{eq:ode_trivial_f_solution}, $F_1({w})$ and $F_2({w})$ are smooth fundamental solutions functions of \cref{eq:ode_trivial_f} (i.e. $F_1({w}_0)=1$, $F_1'({w}_0)=0$, $F_2({w}_0)=0$, and $F_2'({w}_0)=1$) and ${a}_i={a}_i({x},{y})$ and ${b}_i={b}_i({x},{y})$, for $i=1,2$, so ${A}={a}_1\vec{{x}}+{a}_2\vec{{y}}$ and ${B}={b}_1\vec{{x}}+{b}_2\vec{{y}}$ are smooth vector fields on $\RR^2$ equipped with coordinates $({x},{y})$.}

{
Using the commutativity of $f_c$ and $\lb{g_c}{f_c}$ we deduce that 
\begin{align*}
    \lb{F_1{A}+F_2{B}}{F_1'{A}+F_2'{B}}
    =\left(F_1F_2'-F_1'F_2\right)\lb{{A}}{{B}}=0.
\end{align*}
\noindent
By $F_1F_2'-F_1'F_2\neq0$ (since $f_c\wedge \lb{g_c}{f_c}\neq0$), we conclude that $\lb{{A}}{{B}}=0$ and{, therefore, there exists a local diffeomorphism $\psi(x,y)$ that simultaneously rectify $A$ and $B$ (seen as vector fields on $\RR^2$). For simplicity, we still denote the new coordinates by $(x,y)$, i.e. we have $\psi_*A=\vec{x}$ and $\psi_*B=\vec{y}$.} {In coordinates $(x,y,w)$,} {the vector fields $(f_c,g_c)$ take the form}
\begin{align*}
{f}_c=F_1(w)\vec{x}+F_2(w)\vec{y} \quad\textrm{and}\quad {g}_c=\vec{w}
\end{align*}
\noindent
{and therefore we conclude that the system $\Sigma_c=(f_c,g_c)$ is trivial.}
} 
\end{proof}
Remark that in our proof {we start with a canonical pair $(f_c,g_c)$ and we render it trivial by constructing a suitable local coordinate system.}
\begin{remark}\label{rk:error_serres_prf}
    The previous theorem was first discovered by Serres in \cite{serres2006Geometryfeedbackclassification}. In the proof of \cite[Theorem 4.3.3]{serres2006Geometryfeedbackclassification} (but also in \cite[Theorem 4.3]{serres2007curvaturefeedbackclassification} and in \cite[Theorem 3.4]{serres2009Controlsystemszero}){, he shows, using his notation, that $\alpha_2=a_2(u,q_2)-q_1$ and $\frac{\partial a_2}{\partial q_2}=b(u)$ and then consider the case $\alpha_2=a_2(u)-q_1$ and not the general case $\alpha_2=b(u)q_2+a_2(u)-q_1$. The proof of \cite[Theorem 4.3.3]{serres2006Geometryfeedbackclassification}, given for the case $b\equiv0$ (which, using our notation, is equivalent to $\nu\equiv0$){,} still provides an inspiring intuition to treat the general case, as done in our proof.}
\end{remark}
{In the following proposition, we {express the structure functions of a trivial system $(T)$ and} give two {canonical} forms of control-affine system that are trivialisable. Both {canonical} forms are expressed using the canonical pair but in different coordinate systems. {For two smooth scalar functions $F(w)$ and $G(w)$, we define their Wronskian as $\mathtt{W}(F,G)=F'G-FG'$. Recall that for any control-affine system $\Sigma=(f,g)$ satisfying \cref{assumption:1} and \cref{assumption:2} we defined, via \cref{eq:structure_functions_definition}, structure functions $k_1,k_2,k_3$ and $\lambda_1,\lambda_2,\lambda_3$.}
\begin{proposition}\label{prop:normal-forms-trivialisable}
    Consider a control-affine system $\Sigma=(f,g)$ and suppose that it satisfies conditions  \cref{eq:conditions_trivialisable} of \cref{thm:characterisation_trivialisable}. {Then, locally, the following hold
    \begin{enumerate}[label=(\roman*),ref=\textit{(\roman*)}]
        \item\label{prop:normal-forms-trivialisable:1} $\Sigma$ admits the normal form $(T)$, that is,
        \begin{align*}
            \Sigma^T\,:\, \left\{\begin{array}{rl}
            \dot{x} &= F_1(w) \\
            \dot{y} &= F_2(w) \\
            \dot{w} &= u
        \end{array}\right. 
        \end{align*}
        \noindent
        whose structure functions are $k_1=k_2=k_3=0$ and $\lambda_1=-\frac{\mathtt{W}(F_1',F_2')}{\mathtt{W}(F_1,F_2)},\,\lambda_2=0,\,\lambda_3=\frac{\mathtt{W}'(F_1,F_2)}{\mathtt{W}(F_1,F_2)}$.
        \item\label{prop:normal-forms-trivialisable:2} $\Sigma$ admits the canonical forms $\Sigma_c^{T,1}$ and $\Sigma_c^{T,2}$ given, respectively, by 
        \begin{align*}
        \Sigma_c^{T,1}\,:\,\left\{\begin{array}{rl}
            \dot{x} &= F_{c,1}(w) \\
            \dot{y} &= F_{c,2}(w) \\
            \dot{w} &= u
        \end{array}\right.\quad \textrm{and}\quad \Sigma_c^{T,2}\,:\,\left\{\begin{array}{rl}
            \dot{x} &= 1+\varepsilon y u \\
            \dot{y} &=\left(x-\nu(w)y\right)u \\
            \dot{w} &= u
        \end{array}\right.,
        \end{align*}
        \noindent
        where $\frac{\mathtt{W}(F_{c,1}',F_{c,2}')}{\mathtt{W}(F_{c,1},F_{c,2})}\equiv \pm1$ and whose invariants are $(\varepsilon^1,\kappa^1,\nu^1)=\left(-\frac{\mathtt{W}(F_{c,1}',F_{c,2}')}{\mathtt{W}(F_{c,1},F_{c,2})},0,\frac{\mathtt{W}'(F_{c,1},F_{c,2})}{\mathtt{W}(F_{c,1},F_{c,2})}\right)$ and $(\varepsilon^2,\kappa^2,\nu^2)=(\varepsilon, 0, \nu(w))$, respectively.
    \end{enumerate}
    }
\end{proposition}
{
\begin{remark}
    Neither the structure functions $k_i$ nor $\lambda_i$ are feedback invariant. Item \cref{prop:normal-forms-trivialisable:1} asserts that for the normal form $\Sigma^T=(f,g)$, all $k_i=0$, so the pair $(f,g)$ is semi-canonical (thus, actually, $f=f_s$) but, in general, it is not canonical since $\lambda_1$ is a non trivial function. Item \cref{prop:normal-forms-trivialisable:2} assures that, given $\Sigma^T=(f,g)$, we can always choose $w_c$-coordinate, as $w=\phi(w_c)$, such that $F_{c,i}=\phi^*F_i$ satisfy $\frac{\mathtt{W}(F_{c,1}',F_{c,2}')}{\mathtt{W}(F_{c,1},F_{c,2})}=\pm1$ and the corresponding pair $(f_c,g_c)$, where $f_c=F_{c,1}\vec{x}+F_{c,2}\vec{y}$ and $g_c=\beta g$, with $\beta=\phi'$, is canonical.
\end{remark}
}
\begin{proof}
    {The normal form presented in item \cref{prop:normal-forms-trivialisable:1} is a direct consequence of \cref{thm:charac_trivial_sys} and it is a straightforward computation to derive the expressions of the structure functions. }{To obtain the canonical form $\Sigma_c^{T,1}$ of item \cref{prop:normal-forms-trivialisable:2}, we consider $\Sigma^T$ and define $g_c=\beta g$, where $\beta=\left|\lambda_1\right|^{-1/2}$, see \cref{prop:existence_pairs}. We choose $w=\phi(\hat{w})$ satisfying $(\phi^{-1})'\beta=1$. Then in the coordinates $(x,y,\hat{w})$, the system $\Sigma^T$ takes the form $\Sigma_c^{T,1}$, where $F_{c,i}=\phi^*F_i$ and whose third  equation reads $\dot{\hat{w}}=\hat{u}$.} 

    {Finally, the} {canonical} form $\Sigma^{T,2}_c$ is a special case of item \cref{prop:normal_forms:1} of \cref{prop:normal_forms} presented in the next section. 
\end{proof}
The two presented {canonical} forms are somehow \emph{dual} to each other. Indeed, both are given in terms of the canonical pair {$(f_c,g_c)$} of $\Sigma$ and for $\Sigma_c^{T,1}$ we adopt coordinates for which the vector field $g^1_c$ is rectified, whereas for $\Sigma_c^{T,2}$ the coordinates are chosen so that $f^2_c$ is rectified. The two normal forms carry complementary informations about the control-affine system $\Sigma$. {The canonical form $\Sigma_c^{T,1}$ exhibits} the trivial nature of $\Sigma$, but its invariants $\varepsilon$ and $\nu$ are not {immediately} visible, and the {canonical} form {$\Sigma_c^{T,2}$} explicitly identifies the invariants $\varepsilon=\pm1$ and $\nu$ but hides the triviality of the system. The two {canonical} forms show that trivial systems depend on a smooth function of one variable: for $\Sigma_c^{T,2}$ it is clearly $\nu(w)$ and for $\Sigma_c^{T,1}$ it is {the function $F_{c,2}(w)$} that determines {$F_{c,1}(w)$} (or, equivalently, the other way around)  through the ODE {$\frac{\mathtt{W}(F_{c,1}',F_{c,2}')}{\mathtt{W}(F_{c,1},F_{c,2})}=\pm1$}.
}
\section{Normal forms of flat and centro-flat control-affine systems on 3D-manifolds} \label{sec:normal_forms}
We have shown that the curvature $\kappa$ and the centro-affine curvature $\nu$ are two functional feedback {equivariants} of control-affine systems{, hence, their properties define non-equivalent classes of systems.} In {this section}, we propose {a normal form for each class of control-systems that is presented in \cref{table:normal_form_flat} below.} {The presented classes describe {all} the cases for which the curvature $\kappa$ and the centro-affine curvature $\nu$ satisfy $\kappa\nu\equiv0$ {together with the particular sub-cases for which, additionally, either $\kappa$ or $\nu$ is constant.}}

\begin{table}[H]
    \centering
    \setlength{\tabcolsep}{1em} 
    \renewcommand{\arraystretch}{1.2}
    \begin{tabularx}{\textwidth}{c|c|X}
        Notation & Name & Properties  \\
        \hline
        $\Sigma^{\varepsilon,\kappa=0{,\nu}}$ & Flat & Curvature $\kappa$ vanishes \\
        {$\Sigma^{\varepsilon,\kappa,\nu=0}$} & {Centro-flat} & {Centro-affine curvature $\nu$ vanishes} \\
        {$\Sigma^{\varepsilon,\kappa'=0,{\nu=0}}$} & {Flat-constant} & {Curvature $\kappa$ is constant and the centro-affine curvature $\nu$ vanishes} \\
        $\Sigma^{\varepsilon,{\kappa=0},{\nu'=0}}$ & {Centro-}flat-constant & Curvature $\kappa$ vanishes, and {the centro-affine curvature} $\nu$ is constant \\
        $\Sigma^{\varepsilon,{\kappa=0,\nu=0}}$ & Completely flat  & Curvatures $\kappa$ and $\nu$ vanish \\
    \end{tabularx}
    \caption{Nomenclature of subclasses of flat and centro-flat control-affine systems}
    \label{table:normal_form_flat}
\end{table}

%
{Each class of control system presented in the above table is denoted by an upper index $I=(a,b,c)$, which is defined as follows. The first element is always $\varepsilon=\pm1$ and emphasises the dependence of the normal forms on the invariant $\varepsilon$; the second element is either $\kappa=0$ to say that the curvature vanishes or $\kappa'=0$ to express that the curvature is constant (this notation is a bit abusive because $\kappa$ is not a function of one variable in general); finally, the third index is either $\nu=0$ or $\nu'=0$ with the same interpretation as previously.}
{The following proposition {provides a normal form $\Sigma_c^I$ for each class of control-affine systems $\Sigma^I$, where the upper multi-index $I$ is one of the five given in \cref{table:normal_form_flat}. The lower index $c$ indicates that all normal forms $\Sigma_c^I$ are expressed using their canonical pairs.} Recall that the structure function $\nu$ is unique up to its sign, i.e. changing $g_c\mapsto -g_c$ yields $\nu\mapsto-\nu$, hence in normal forms below we suppose that $\nu\geq0$. }
\begin{proposition}[Normal forms {of flat control-affine systems}]\label{prop:normal_forms}
Consider a control-affine system {$\Sigma=(f,g)$} together with {its} invariants $\varepsilon$, $\kappa$, and $\nu$. Then, the following statements hold locally (all normal forms below are {represented by a canonical pair {$(f_c,g_c)$} and} considered around an arbitrary point $(x_0,y_0,w_0)\in\RR^3$).
\begin{enumerate}[label=(\roman*),ref=\textit{(\roman*)}]
    \item \label{prop:normal_forms:1} If $\kappa=0$, then $\Sigma$ is  {locally} feedback equivalent to
    \begin{align*}
    \Sigma^{\varepsilon,\kappa=0{,\nu}}_c\,:\,\left\{\begin{array}{rll}
            \dot{x} &= 1 +& {a(y,w)u} \\ 
            \dot{y} &=   & {\left(x+b(y,w)\right)u} \\ 
            \dot{w} &=   & {c(y,w)u} 
        \end{array}\right.,
    \end{align*}
    \noindent
    whose invariants are $\varepsilon$, $\kappa=0$, and $\nu=\nu_1(y,w)x+\nu_0(y,w)$, {and the functions $a, b$ and $c$ satisfy the following differential equations 
	\begin{align*}
	\frac{\partial a}{\partial y}&=\varepsilon+\nu_1(y,w)a(y,w),\quad &a(y_0,w)=0,\\
    \frac{\partial b}{\partial y}&=\nu_1(y,w)b(y,w)-\nu_0(y,w),\quad &b(y_0,w)=0,\\
    \frac{\partial c}{\partial y}&=\nu_1(y,w)c(y,w),\quad &c(y_0,w)=1,
	\end{align*}	 
	\noindent
	and thus are given by 
	\begin{align*}
	 a(y,w)&=\left[\varepsilon\int_{y_0}^y\exp\left(-\int_{y_0}^{\tau}\nu_1(t,w)\diff t\right)\diff \tau\right]\exp\left(\int_{y_0}^y\nu_1(\tau,w)\diff \tau\right),\\ 
	 b(y,w)&=\left[-\int_{y_0}^y\nu_0(\tau,w)\exp\left(-\int_{y_0}^{\tau}\nu_1(t,w)\diff t\right)\diff\tau\right]\exp\left(\int_{y_0}^y\nu_1(\tau,w)\diff \tau\right),\\
        c(y,w)&=\exp\left(\int_{y_0}^y\nu_1(\tau,w)\diff \tau\right).
	\end{align*}	   
    }
    \item \label{prop:normal_forms:2} {If $\nu =0$, then $\Sigma$ is {locally} feedback equivalent to 
    \begin{align*}
        \Sigma^{\varepsilon{,\kappa},\nu=0}_c\,:\,\left\{\begin{array}{rl}
            \dot{x} &= r(x,y)\,c_{\varepsilon}(w) \\ 
            \dot{y} &= r(x,y)\,s_{\varepsilon}(w) \\ 
            \dot{w} &= \varepsilon \frac{\partial r}{\partial y}\,c_{\varepsilon}(w)+\frac{\partial r}{\partial x}\,s_{\varepsilon}(w) + u \\
        \end{array}\right., 
    \end{align*}
    \noindent
    where 
    \begin{align*}
        c_{\varepsilon}(w) = \frac{e^{w\sqrt{\varepsilon}} + e^{-w\sqrt{\varepsilon}} }{2} \quad \textrm{and} \quad s_{\varepsilon}(w)=\frac{e^{w\sqrt{\varepsilon}} -e^{-w\sqrt{\varepsilon}} }{2\sqrt{\varepsilon}},
    \end{align*}
    \noindent
    whose invariants are $\varepsilon$, $\kappa=\kappa(x,y)$, $\nu=0$, and the function $r(x,y)$ satisfies $r>0$ and the following non-linear partial differential equation
    \begin{align}\label{eq:relation:r_curvature_kappa}
       -r(x,y)^2\left( \frac{\partial^2}{\partial x^2}-\varepsilon\frac{\partial^2}{\partial y^2}\right)(\ln r(x,y))= \kappa(x,y).
    \end{align}
    }
    \item \label{prop:normal_forms:3} {If $\kappa$ and $\nu$ are constant then 
    \begin{align}\label{eq:relation:curvatures:constants}
        \kappa\nu=0,
    \end{align}
    \noindent
    i.e. at least one of them {vanishes}.}
    \item \label{prop:normal_forms:4} {{If $\kappa=0$ and $\nu$ is constant}, then $\Sigma$ is locally feedback equivalent to} 
    \begin{enumerate}[label=(\alph*),ref=\textit{(\alph*)}]
        \item\label{prop:normal_forms:4:a} If $\varepsilon=1$, then
        \begin{align*}
            \Sigma_c^{+,{\kappa=}0,{\nu'=0}}\,:\,\left\{
            \begin{array}{rl}
                \dot{x} &= e^{\nu w}e^{w\sqrt{\nu^2+4}}\\
                \dot{y} &= e^{\nu w}e^{-w\sqrt{\nu^2+4}}  \\
                \dot{w} &= \frac{1}{2} u
            \end{array}
            \right.,\quad\textrm{where}\quad \nu\geq0;
        \end{align*}
        \item\label{prop:normal_forms:4:b} If $\varepsilon=-1$ and $\nu>2$, then 
        \begin{align*}
            \Sigma_c^{-,{\kappa=}0,{\nu'=0},+}\,:\,\left\{
            \begin{array}{rl}
                \dot{x} &= e^{\nu w}e^{w\sqrt{\nu^2-4}}\\
                \dot{y} &= e^{\nu w}e^{-w\sqrt{\nu^2-4}}  \\
                \dot{w} &= \frac{1}{2}u
            \end{array}
            \right.;
        \end{align*}
        \item \label{prop:normal_forms:4:c}If $\varepsilon=-1$ and $\nu=2$, then
        \begin{align*}
            \Sigma_c^{-,{\kappa=}0,{\nu'=0},0}\,:\,\left\{
            \begin{array}{rl}
                \dot{x} &= e^w \\
                \dot{y} &= we^w \\
                \dot{w} &= u
            \end{array}
            \right.;
        \end{align*}
        \noindent
        \item \label{prop:normal_forms:4:d}If $\varepsilon=-1$ and $0\leq\nu<2$, then 
        \begin{align*}
            \Sigma_c^{-,{\kappa=}0,{\nu'=0},-}\,:\,\left\{
            \begin{array}{rl}
                \dot{x} &= e^{\nu w}\cos\left({w\sqrt{4-\nu^2}}\right)\\
                \dot{y} &= e^{\nu w}\sin\left({w\sqrt{4-\nu^2}}\right) \\
                \dot{w} &= \frac{1}{2}u
            \end{array}
            \right.
        \end{align*}
    \end{enumerate}
    Moreover, for the four normal forms above, $\varepsilon$, $\kappa=0$, and $\nu$ are complete invariants.    
    \item \label{prop:normal_forms:5} {{If $\nu=0$ and $\kappa$ is constant}, then $\Sigma$ is locally feedback equivalent to
    \begin{align*}
         \Sigma^{\varepsilon,\kappa'=0,\nu=0}_c\,:\,\left\{\begin{array}{rl}
            \dot{x} &= \left(1-\frac{\kappa}{4}\left(x^2-\varepsilon y^2\right)\right)\,c_{\varepsilon}(w) \\ 
            \dot{y} &= \left(1-\frac{\kappa}{4}\left(x^2-\varepsilon y^2\right)\right)\,s_{\varepsilon}(w) \\
            \dot{w} &= \frac{-\kappa}{2} \left(y\,c_{\varepsilon}(w)-x\,s_{\varepsilon}(w)\right) + u \\ 
        \end{array}\right.,
    \end{align*}
    whose complete invariants are $\varepsilon$, $\kappa$, and $\nu=0$.
    }
    \item \label{prop:normal_forms:6} If $\kappa=0$ and $\nu=0$, then $\Sigma$ is  {locally} feedback equivalent to 
    \begin{align*}
        \Sigma^{\varepsilon,{\kappa=}0,{\nu=}0}_c\,:\,\left\{\begin{array}{rl}
            \dot{x} &= c_{\varepsilon}(w)\\
            \dot{y} &= s_{\varepsilon}(w)\\
            \dot{w} &= u 
        \end{array}\right..
    \end{align*}
\end{enumerate}
\end{proposition}
{Before {presenting} a proof for those normal forms, we give some remarks about them. For item \cref{prop:normal_forms:1} we adopt coordinates, where the vector field $f_c$ is rectified, whereas for the {other} normal forms we choose coordinates in which the vector field $g_c$ is rectified.} 
{The first normal form $\Sigma_c^{\varepsilon,\kappa=0, \nu}$ of flat control-affine systems describes the most general form of a system for which the curvature $\kappa$ vanishes. On the other hand, the normal form $\Sigma_c^{\varepsilon,\kappa,\nu=0}$ of item \cref{prop:normal_forms:2}, describes systems for which the centro-affine curvature $\nu$ vanishes. All other items are then special cases of those two general normal forms.} 

Recall that $\nu$ is unique up to its sign, that is why {in item \cref{prop:normal_forms:4}} we have $\nu\geq0$ for \cref{prop:normal_forms:4}-\cref{prop:normal_forms:4:a} to \cref{prop:normal_forms:4}-\cref{prop:normal_forms:4:d}. {It is remarkable that if $\kappa$ and $\nu$ are constant (hence true invariants) then at least one of them is zero as asserted in item \cref{prop:normal_forms:3}. Moreover, relation \cref{eq:relation:curvatures:constants} already appeared in \cite{wilkens1998Centroaffinegeometryplane}, where the four families of normal forms given by $\kappa=0$ and $\nu$ constant were listed (but the non-invariance of the sign of $\nu$ was not discussed there).} 
The two normal forms {of} \cref{prop:normal_forms:6} with {$\varepsilon=\pm1$ and } $\kappa=\nu=0$ are , respectively, given by
\begin{align*}
    \Sigma^{+,\kappa=0,\nu=0}_{c}\,:\,\left\{\begin{array}{rl}
            \dot{x} &=\cosh(w) \\
            \dot{y} &=\sinh(w)  \\
            \dot{w} &= u  
        \end{array}\right.\quad\textrm{and}\quad
        \Sigma^{-,\kappa=0,\nu=0}_{c}\,:\,\left\{\begin{array}{rl}
            \dot{x} &=\cos(w) \\
            \dot{y} &=\sin(w)  \\
            \dot{w} &= u  
        \end{array}\right.
\end{align*}
\noindent
correspond to hyperbolic and elliptic systems without parameters and have been extensively analysed and differently characterised in \cite{schmoderer2021Conicnonholonomicconstraints,SCHMODERER2022105397}.
\begin{proof}\leavevmode
{For each item, we consider a control-affine system $\Sigma_c=(f_c,g_c)$ given by the canonical pair and with $\varepsilon=\pm1$ and structure functions $\kappa$ and $\nu$.}
\begin{enumerate}[label=\textit{(\roman*)},ref=\textit{(\roman*)}]
    \item \label{prf:normal_forms:1} 
     Since $\kappa=0$, using relation \cref{eq:canonical_structure_constants}, we conclude that the vector fields $f_{{c}}$ and $\lb{f_{{c}}}{g_{{c}}}$ are commuting. Therefore, we can rectify them simultaneously {to get} $f_{{c}}=\vec{x}$ and $\lb{f_c}{g_{{c}}}=\vec{y}$. Afterwards, we determine the form of the vector field $g_{{c}}$. First, it satisfies $\lb{\vec{x}}{g_{{c}}}=\vec{y}$ and thus we immediately conclude 
    \begin{align*}
        g_{{c}}=a(y,w)\vec{x}+(x+b(y,w))\vec{y}+c(y,w)\vec{w}.
    \end{align*}
    \noindent
    Moreover, assumptions \cref{assumption:1} and \cref{assumption:2} imply that $c\neq0$ and $\frac{\partial }{\partial y}\left(\frac{a}{c}\right)\neq0$. Second, {for} $g_c$ {we have} $\lb{g_{{c}}}{-\vec{y}}=\varepsilon \vec{x}+\mu g_c-\nu\vec{y}$, where {the functions} $\mu$ and $\nu$ {satisfy} \cref{eq:relation_struct_functions:2_p} and \cref{eq:relation_struct_functions:3_p} {and therefore $\nu=\nu_1(y,w)x+\nu_0(y,w)$ and $\mu=\nu_1(y,w)$}. Hence, the functions $a$, $b$, and $c$ of $g_{{c}}$ satisfy 
    \begin{align*}
        \frac{\partial a}{\partial y}(y,w)&=\varepsilon+\nu_1(y,w)a(y,w),\\
        \frac{\partial b}{\partial y}(y,w)&=\nu_1(y,w)b(y,w)-\nu_0(y,w),\\
        \frac{\partial c}{\partial y}(y,w)&=\nu_1(y,w)c(y,w).
    \end{align*}
    \noindent
    Solutions of those equations are, respectively,
    \begin{subequations}
    \begin{align}
        a(y,w)&=\left[\varepsilon\int_{y_0}^y\exp\left(-\int_{y_0}^{\tau}\nu_1(t,w)\diff t\right)\diff \tau + A(w)\right]\exp\left(\int_{y_0}^y\nu_1(\tau,w)\diff \tau\right), \label{eq:sol:normal_form_flat_canonic:a}\\ 
        b(y,w)&=\left[-\int_{y_0}^y\nu_0(\tau,w)\exp\left(-\int_{y_0}^{\tau}\nu_1(t,w)\diff t\right)\diff\tau+B(w)\right]\exp\left(\int_{y_0}^y\nu_1(\tau,w)\diff \tau\right),\label{eq:sol:normal_form_flat_canonic:b}\\
        c(y,w)&=C(w)\exp\left(\int_{y_0}^y\nu_1(\tau,w)\diff \tau\right).\label{eq:sol:normal_form_flat_canonic:c}
    \end{align}
    \end{subequations}
    \noindent
    Updating the coordinates, we can set $C(w)=1$, $A(w)=B(w)=0$; and in those coordinates we obtain the normal form $\Sigma_c^{\varepsilon,\kappa=0,\nu}$. 
    \item {Suppose that $\nu=0$ and choose coordinates $(\bar{x},\bar{y},\bar{w})$ such that $g_c=\vec{\bar{w}}$. Then, by relation \cref{eq:canonical_structure_constants} we conclude that 
    \begin{align*}
        f_c = \bar{A}(\bar{x},\bar{y})c_{\varepsilon}(\bar{w})+\bar{B}(\bar{x},\bar{y})s_{\varepsilon}(\bar{w}), \quad c_{\varepsilon}(\bar{w}) = \frac{e^{\bar{w}\sqrt{\varepsilon}}+e^{-\bar{w}\sqrt{\varepsilon}}}{2}, \quad s_{\varepsilon}(\bar{w}) = \frac{e^{\bar{w}\sqrt{\varepsilon}}-e^{-\bar{w}\sqrt{\varepsilon}}}{2\sqrt{\varepsilon}},
    \end{align*}
    \noindent
    where $\bar{A}=a_1\Vec{\bar{x}}+a_2\Vec{\bar{y}}+a_3\Vec{\bar{w}}$, with $a_i=a_i(\bar{x},\bar{y})$, and $\bar{B}=b_1\Vec{\bar{x}}+b_2\Vec{\bar{y}}+b_3\Vec{\bar{w}}$, with $b_i=b_i(\bar{x},\bar{y})$ are smooth vector fields. By assumption \cref{assumption:1}, we conclude that $a_1b_2-a_2b_1\neq0$, hence $\bar{\bar{A}}=a_1\Vec{\bar{x}}+a_2\Vec{\bar{y}}$ and $\bar{\bar{B}}=b_1\Vec{\bar{x}}+b_2\Vec{\bar{y}}$ form a moving frame of the tangent bundle of $\XXX=\mathcal{O}/\GGG$, where $\mathcal{O}$ is an open subset of $\RR^3$, in which the rectifying coordinates $(\bar{x},\bar{y},\bar{w})$ are defined, and $\GGG=\distrib{\vec{\bar{w}}}$. We define a metric $\bar{\bar{\textsf{g}}}$ on $\XXX$ by declaring $(\bar{\bar{A}},\bar{\bar{B}})$ orthonormal, i.e. 
    \begin{align*}
        \bar{\bar{\textsf{g}}}(\bar{\bar{A}},\bar{\bar{A}})=1, \quad \bar{\bar{\textsf{g}}}(\bar{\bar{A}},\bar{\bar{B}})=0, \quad\textrm{and} \quad\bar{\bar{\textsf{g}}}(\bar{\bar{B}},\bar{\bar{B}})=-\varepsilon. 
    \end{align*}
    \noindent
    Notice that the signature of $\bar{\bar{\textsf{g}}}$ is $(+,-\sgn{\varepsilon})$, hence $\bar{\bar{\textsf{g}}}$ is definite for $\varepsilon=-1$ and indefinite for $\varepsilon=1$. Since all metrics on $2$-dimensional manifolds are locally conformally flat, we conclude that there exists an isometry $(x,y)=\psi(\bar{x},\bar{y})$ such that $\bar{\bar{\textsf{g}}}=\psi^*\textsf{g}$, where $\textsf{g}=\varrho(x,y)(\diff x^2-\varepsilon \diff y^2)$ with $\varrho>0$. In {the} coordinates $(x,y)$, both {the pair} $(\tilde{A},\tilde{B})$, with $\Tilde{A}=\psi_*\bar{\bar{A}}$ and $\Tilde{B}=\psi_*\bar{\bar{B}}$, and $\left(\vec{x},\vec{y}\right)$ form an orthornormal frame for $\textsf{g}$ so we have 
    \begin{align*}
        (\tilde{A},\tilde{B}) = r(x,y) I(x,y) \left(\vec{x},\vec{y}\right),
    \end{align*}
    \noindent
    where $r=\frac{1}{\sqrt{\varrho}}$ and $I(x,y)$ is a {a linear} isometry, i.e. it belongs to the (pseudo)-orthonormal group $O(1,-\varepsilon)$. Using a {suitable} change of {the} variable ${w}=\bar{w}+h({x},{y})$ we can {get rid of} $I(x,y)$. Finally, in coordinates $(x,y,w)$, the vector field $f_c$ of the control system takes the form 
    \begin{align*}
        f_c=r(x,y)c_{\varepsilon}(w)\Vec{x}+r(x,y)s_{\varepsilon}(w)\Vec{y} + \left(a(x,y)c_{\varepsilon}(w)+b(x,y)s_{\varepsilon}(w)\right)\Vec{w}.
    \end{align*}
    \noindent
    We now use the structure equations \cref{eq:canonical_structure_constants} and deduce that necessarily
    \begin{align*}
        a(x,y)= \varepsilon \frac{\partial r}{\partial y}\quad \textrm{and}\quad b(x,y)=  \frac{\partial r}{\partial x}
    \end{align*}
    \noindent
    and that $r$ satisfies 
    \begin{align*}
        \varepsilon r\frac{\partial ^2r}{\partial y^2} - \varepsilon \left(\frac{\partial r}{\partial y}\right)^2+\left(\frac{\partial r}{\partial x}\right)^2-r\frac{\partial ^2r}{\partial x^2} = \kappa (x,y). 
    \end{align*}
    \noindent
    which can be expressed in the form of \cref{eq:relation:r_curvature_kappa}.
    }
    \item {If $\kappa$ and $\nu$ are constants, then due to relation \cref{eq:relation_kappa_nu}, we immediately conclude \cref{eq:relation:curvatures:constants}. }
    \item \label{prf:normal_form:constant}
    Suppose that $\kappa=0$ and $\nu$ is constant, then $\Sigma$ satisfies condition \cref{eq:conditions_trivialisable} of \cref{thm:characterisation_trivialisable} and thus $\Sigma$ is locally trivialisable. {Using the results of item \cref{prop:normal-forms-trivialisable:2} of \cref{prop:normal-forms-trivialisable}, we take $\Sigma$ in the form of $\Sigma_c^{T,1}$ for which $f_c=F_{c,1}(w)\vec{x}+F_{c,2}(w)\vec{y}$ and $g_c=\vec{w}$ form a canonical pair.} Using \cref{eq:canonical_structure_constants}, we conclude that the functions $F_{c,i}$, for $i=1,2$, satisfy the following second order linear ordinary differential equation
    \begin{align}\label{eq:prf:normal_forms_trivial}
        F_{c,i}''(w)=\varepsilon F_{c,i}(w)+\nu F_{c,i}'(w). 
    \end{align}
    \noindent
    Solutions are dictated by the sign of the discriminant $\Delta=\nu^2+4\varepsilon$ of the characteristic polynomial of the ODE. Moreover, the roots of the characteristic polynomial are $r_{1/2}=\frac{\nu\pm\sqrt{\Delta}}{2}$. Recall that the sign of $\nu$ is not invariant and thus by choosing $w$ suitably we can always get $\nu\geq0$. Moreover, it is a trivial calculation to check that the solutions given below are fundamental solutions of \cref{eq:prf:normal_forms_trivial}, i.e. we just need to compute the Wronskian at $w_0$. 
    \begin{enumerate}[label=(\alph*),ref=(\alph*)]
        \item \label{prf:normal_form:constant:a} If $\varepsilon=+1$, then $\Delta>0$ for all $\nu\geq0$. Solutions of \cref{eq:prf:normal_forms_trivial} are given by (after normalising $w$ with $\frac{1}{2}$)
        \begin{align*}
            F_{c,1}(w)=e^{\nu w}e^{w\sqrt{\nu^2+4}}\quad\textrm{and}\quad F_{c,2}(w)=e^{\nu w}e^{-w\sqrt{\nu^2+4}},
        \end{align*}
        \noindent
        and we obtain the normal form $\Sigma_c^{+,\kappa=0,\nu'=0}$.
        \item If $\varepsilon=-1$ and $\nu>2$, then $\Delta>0$, and solving \cref{eq:prf:normal_forms_trivial} gives 
        \begin{align*}
            F_{c,1}(w)=e^{\nu w}e^{w\sqrt{\nu^2-4}},\quad\textrm{and}\quad F_{c,2}(w)=e^{\nu w}e^{-w\sqrt{\nu^2-4}},
        \end{align*}
        \noindent
        and we obtain the normal form $\Sigma_c^{-,{\kappa=}0,{\nu'=0},+}$.
        \item If $\varepsilon=-1$ and $\nu=2$, then $\Delta=0$, and the solutions of \cref{eq:prf:normal_forms_trivial} are 
        \begin{align*}
            F_{c,1}(w)=w e^{w},\quad\textrm{and}\quad F_{c,2}(w)=e^{ w},
        \end{align*}
        \noindent
        which gives $\Sigma_c^{-,{\kappa=}0,{\nu'=0},0}$.
        \item\label{prf:normal_form:constant:d} If $\varepsilon=-1$ and $0\leq\nu<2$, then $\Delta<0$, and the solutions of \cref{eq:prf:normal_forms_trivial} are 
        \begin{align*}
            F_{c,1}(w)=e^{\nu w}\cos\left(w\sqrt{4-\nu^2}\right),\quad\textrm{and}\quad F_{c,2}(w)=e^{\nu w}\sin\left(w\sqrt{4-\nu^2}\right),
        \end{align*}
        \noindent
        which gives $\Sigma_c^{-,{\kappa=}0,{\nu'=0},-}$.
    \end{enumerate}
    \item \label{prf:normal_form:constant_kappa} {Assume that $\kappa$ is constant and $\nu=0$, then we refine the normal form $\Sigma_c^{\varepsilon,\kappa,\nu=0}$ of item \cref{prop:normal_forms:2}. We recognize that equation \cref{eq:relation:r_curvature_kappa} satisfied by $r(x,y)$ describes the curvature (in the usual differential geometry sense) of the metric $\textsf{g}=\frac{1}{r^2}(\diff x^2-\varepsilon \diff y^2)$. By assumption, the curvature of $\textsf{g}$ is constant (equal to $-\kappa$)  and by Minding's theorem, surfaces with the same constant curvature are locally isometric. Therefore, there exists an isometry $(\tilde{x},\tilde{y})=\psi(x,y)$ such that $\textsf{g}=\psi^*\tilde{\textsf{g}}$ with 
    \begin{align*}
      \tilde{\textsf{g}} =  \left( \frac{1}{1-\frac{\kappa}{4}\left(\tilde{x}^2-\varepsilon \tilde{y}^2\right)}\right)^2\left(\diff \tilde{x}^2-\varepsilon \diff \tilde{y}^2\right),
    \end{align*}
    \noindent
    which is also of curvature $-\kappa$. The action of the isometry on $(\dot{x},\dot{y})$ can be compensated by applying $w\mapsto w+h(x,y)$, for a suitable function $h$, thus we obtain that the system takes the form of $\Sigma_c^{\varepsilon,\kappa,\nu=0}$ with $r(x,y)=1-\frac{\kappa}{4}\left(x^2-\varepsilon y^2\right)$, i.e. we get $\Sigma_c^{\varepsilon,\kappa'=0,\nu=0}$.
    }
    \item The normal forms $\Sigma_c^{\varepsilon,\kappa=0,\nu=0}$ is a special case of item \cref{prf:normal_form:constant_kappa} with $\kappa=0$. 
\end{enumerate}
\end{proof}

\section{Conclusions and Perspectives}

In this paper, we have analysed in details the notion of triviality adapted to the context of control-affine systems. We proposed two new characterisations of trivial control-affine system, one of them is based on the existence of an abelian subalgebra of the Lie algebra of infinitesimal symmetries. In particular, we gave a normal form of trivial control-affine systems for which the Lie algebra of infinitesimal symmetries has a transitive almost abelian Lie subalgebra. In the future, we will be interested in extending our result to the case of multi-input systems and we will try to propose other characterisation of triviality that are purely geometric. In the second part of the paper, we have revisited results due to Serres \cite{serres2006Geometryfeedbackclassification} and we give novel proof of his characterisation of trivial systems on $3$-dimensional manifolds with scalar inputs. In particular, our characterisation uses a discrete invariant $\varepsilon=\pm1$ and two well-defined functional invariants of {feedback transformations}: the curvature $\kappa$ {(introduced by Agrachev \cite{agrachev1997Feedbackinvariantoptimalcontrol})} and the centro-affine curvature $\nu$ (studied by Wilkens \cite{wilkens1998Centroaffinegeometryplane}). We show that those invariants can explicitly be computed for any control-affine system and that a canonical pair of vector fields $(f_c,g_c)$, on which $\kappa$ and $\nu$ appear explicitly, can also be constructed with a purely algebraically defined feedback transformation. Then, we extended the results of Serres and Wilkens by giving several normal forms of control-affine systems. In the future, our goal is two-folds: first we will be interested in the question {of how to enlarge the triple $(\varepsilon, \kappa,\nu)$ to a set of} \emph{complete invariants} of control-affine systems (on $3$D manifolds with scalar control). {Identifying a set of complete invariants would be helpful in obtaining normal forms of control-affine system in dimension three}. Second, we will be interested in generalising our characterisation of trivial control-affine systems to the multi-input case, in particular the {notion of curvature of dynamics pairs, as proposed in \cite{jakubczyk2013Vectorfieldsdistributions},} seems promising. 

\setcounter{biburllcpenalty}{7000}
\setcounter{biburlucpenalty}{8000}
\printbibliography

\appendix
\section{Detailed computations for \texorpdfstring{\cref{lem:transform_structure_functions}}{Lemma \ref*{lem:transform_structure_functions}}}\label{apdx:transform_structure_functions}
In this appendix, we detail the computation to obtain relations \cref{eq:relation_struct_functions:1}-\cref{eq:relation_struct_functions:3} between structure functions and we prove transformation rules \cref{eq:transform_struct_funct_k} and \cref{eq:transform_struct_funct_lambda} that show how the structure functions are changed under a feedback transformation. Consider a control-affine system $\Sigma=(f,g)$ with structure functions $(k_1,k_2,k_3)$ and $(\lambda_1,\lambda_2,\lambda_3)$.

First, by applying the Jacobi identity to $\lb{f}{\lb{g}{\lb{g}{f}}}$ we deduce that $\lb{f}{\lb{g}{\lb{g}{f}}}=-\lb{g}{\lb{f}{\lb{f}{g}}}$. We compute the left-hand-side and the right-hand-side separately: 
\begin{align*}
    \lb{f}{\lb{g}{\lb{g}{f}}}&=\dL{f}{\lambda_1}f+\dL{f}{\lambda_2}g-\lambda_2\lb{g}{f}+\dL{f}{\lambda_3}-\lambda_3\left(k_1g+k_2\lb{g}{f}+k_3\lb{g}{\lb{g}{f}}\right),\\
    &= \dL{f}{\lambda_1}f+\left(\dL{f}{\lambda_2}-\lambda_3 k_1\right)g+\left(\dL{f}{\lambda_3}-\lambda_2-\lambda_3 k_2\right)\lb{g}{f}\\
    &\quad-\lambda_3 k_3\left(\lambda_1 f+\lambda_2 g+\lambda_3\lb{g}{f}\right),\\ %
    &= \left(\dL{f}{\lambda_1}-\lambda_3 k_3\lambda_1\right)f+\left(\dL{f}{\lambda_2}-\lambda_3 k_1-\lambda_3 k_3\lambda_2\right)g\\
    &\quad+\left(\dL{f}{\lambda_3}-\lambda_2-\lambda_3 k_2-\lambda_3^2 k_3\right)\lb{g}{f}.
\end{align*}
\noindent
And on the other hand we have 
\begin{align*}
    \lb{g}{\lb{f}{\lb{f}{g}}} &= \dL{g}{k_1}g+\dL{g}{k_2}\lb{g}{f}+k_2\left(\lambda_1 f+\lambda_2 g+\lambda_3\lb{g}{f}\right)\\
    &\quad+\dL{g}{k_3}\left(\lambda_1 f+\lambda_2 g+\lambda_3\lb{g}{f}\right) +k_3\lb{g}{\lambda_1 f+\lambda_2 g+\lambda_3\lb{g}{f}},\\ 
    &= \left(k_2\lambda_1+\lambda_1\dL{g}{k_3}\right)f+\left(\dL{g}{k_1}+k_2\lambda_2+\dL{g}{k_3}\lambda_2\right)g\\
    &\quad+\left(\dL{g}{k_2}+k_2\lambda_3+\dL{g}{k_3}\lambda_3\right)\lb{g}{f}\\ 
    & \quad +k_3\left(\dL{g}{\lambda_1}f+\lambda_1\lb{g}{f}+\dL{g}{\lambda_2}g+\dL{g}{\lambda_3}\lb{g}{f}+\lambda_3\left(\lambda_1 f+\lambda_3 g+\lambda_3\lb{g}{f}\right)\right),\\
    &= \left(k_2\lambda_1+\lambda_1\dL{g}{k_3}+k_3\dL{g}{\lambda_1}+k_3\lambda_3\lambda_1\right)f\\
    &\quad +\left(\dL{g}{k_1}+k_2\lambda_2+\dL{g}{k_3}\lambda_2+k_3\dL{g}{\lambda_2}+k_3\lambda_3\lambda_2\right)g\\
    & \quad +\left(\dL{g}{k_2}+k_2\lambda_3+\dL{g}{k_3}\lambda_3+k_3\lambda_1+k_3\dL{g}{\lambda_3}+k_3\lambda_3^2\right)\lb{g}{f}
\end{align*}
\noindent
Identifying the terms in front of $f$, $g$, and $\lb{g}{f}$ we obtain equations \cref{eq:relation_struct_functions:1} to \cref{eq:relation_struct_functions:3}. 

Now, we apply a feedback transformation of the form $\tilde{f}=f+g\alpha$ and $\tilde{g}=g\beta$ and we get first $\lb{\tilde{f}}{\tilde{g}}=\beta\lb{f}{g}+\gamma g$, where $\gamma=\dL{f}{\beta}+\alpha\dL{g}{\beta}-\beta\dL{g}{\alpha}$. Second, 
\begin{align*}
    \lb{\tilde{g}}{\lb{\tilde{g}}{\tilde{f}}} &= \beta^2\lb{g}{\lb{g}{f}}+\beta\dL{g}{\beta}\lb{g}{f}+\left(-\beta\dL{\lb{g}{f}}{\beta}-\beta\dL{g}{\gamma}+\gamma\dL{g}{\beta}\right)g \\ 
    &=\beta^2\lambda_1 f + \left(\beta^2\lambda_2-\beta\dL{\lb{g}{f}}{\beta}-\beta\dL{g}{\gamma}+\gamma\dL{g}{\beta}\right)g +\left(\beta^2\lambda_3+\beta\dL{g}{\beta}\right)\lb{g}{f} \\ 
    &=\beta^2\lambda_1\tilde{f}+\left(-\beta^2\lambda_1\alpha +\beta^2\lambda_2-\beta\dL{\lb{g}{f}}{\beta}-\beta\dL{g}{\gamma}+\gamma\dL{g}{\beta}+\gamma\left(\beta\lambda_3+\dL{g}{\beta}\right)\right)g\\
    &\quad+ \left(\beta\lambda_3+\dL{g}{\beta}\right)\lb{\tilde{g}}{\tilde{f}},
\end{align*}
\noindent
implying that $\tilde{\lambda}_1=\beta^2\lambda_1$, $\tilde{\lambda}_2=\beta\lambda_2-\beta\lambda_1\alpha+\gamma\lambda_3  -\dL{\lb{g}{f}}{\beta}-\dL{g}{\gamma}+2\gamma\dL{g}{\ln|\beta|}$, and $\tilde{\lambda}_3=\beta\lambda_3+\dL{g}{\beta}$. Third, we have 
\begin{align*}
    \lb{\tilde{f}}{\lb{\tilde{f}}{\tilde{g}}} &= \lb{f+g\alpha}{\beta\lb{f}{g}+\gamma g}\\
    &=\beta\lb{f}{\lb{f}{g}}+\dL{f}{\beta}\lb{f}{g}+\dL{f}{\gamma}g+\gamma\lb{f}{g}\\
    &\quad+\alpha\beta\lb{g}{\lb{f}{g}}+\alpha\dL{g}{\beta}\lb{f}{g}-\beta\dL{\lb{f}{g}}{\alpha}g+\alpha\dL{g}{\gamma}g-\gamma\dL{g}{\alpha}g,\\
    &= \left(\beta k_1+\dL{f}{\gamma}+\beta\dL{\lb{g}{f}}{\alpha}+\alpha\dL{g}{\gamma}-\gamma\dL{g}{\alpha}\right)g + \left(\beta k_2-\dL{f}{\beta}-\gamma-\alpha\dL{g}{\beta}\right)\lb{g}{f}\\
    &\quad+\left(\beta k_3-\alpha\beta\right)\lb{g}{\lb{g}{f}},\\ 
    &=\left(\beta k_1+\dL{f}{\gamma}+\beta\dL{\lb{g}{f}}{\alpha}+\alpha\dL{g}{\gamma}-\gamma\dL{g}{\alpha}\right)g + \left(\beta k_2-\dL{f}{\beta}-\gamma-\alpha\dL{g}{\beta}\right)\lb{g}{f}\\
    &\quad+ \frac{1}{\beta^2}\left(\beta k_3-\alpha\beta\right)\left\{\lb{\tilde{g}}{\lb{\tilde{g}}{\tilde{f}}}-\beta\dL{g}{\beta}\lb{g}{f}+\beta\dL{\lb{g}{f}}{\beta}g+\beta\dL{g}{\gamma}g-\gamma\dL{g}{\beta}g\right\},
\end{align*}
\noindent
implying $\tilde{k}_3=\frac{1}{\beta}\left(k_3-\alpha\right)$. Next, continuing the computation (denoting $\tilde{h}=\tilde{k}_3\lb{\tilde{g}}{\lb{\tilde{g}}{\tilde{f}}}$):
\begin{align*}
    \lb{\tilde{f}}{\lb{\tilde{f}}{\tilde{g}}} &= \left(\beta k_1+\dL{f}{\gamma}+\beta\dL{\lb{g}{f}}{\alpha}+\alpha\dL{g}{\gamma}-\gamma\dL{g}{\alpha}+\tilde{k}_3\left(\beta\dL{\lb{g}{f}}{\beta}+\beta\dL{g}{\gamma}-\gamma\dL{g}{\beta}\right)\right)g +\tilde{h} \\
    &\quad + \left(\beta k_2-\dL{f}{\beta}-\gamma-\alpha\dL{g}{\beta}-\tilde{k}_3\beta\dL{g}{\beta}\right)\lb{g}{f}+\tilde{h}  \\
    &= \left(\beta k_1+\dL{f}{\gamma}+\beta\dL{\lb{g}{f}}{\alpha}+\alpha\dL{g}{\gamma}-\gamma\dL{g}{\alpha}+\tilde{k}_3\left(\beta\dL{\lb{g}{f}}{\beta}+\beta\dL{g}{\gamma}-\gamma\dL{g}{\beta}\right)\right)g \\ 
    &\quad \frac{1}{\beta}\left(\beta k_2-\dL{f}{\beta}-\gamma-\alpha\dL{g}{\beta}-\tilde{k}_3\beta\dL{g}{\beta}\right)\left\{\lb{\tilde{g}}{\tilde{f}}+\gamma g\right\}+\tilde{h} ,
\end{align*}
\noindent
implying $\tilde{k}_2=k_2-\dL{f}{\ln|\beta|}-\frac{\gamma}{\beta}-\alpha\dL{g}{\ln|\beta|}-\tilde{k}_3\dL{g}{\beta}$ and finally
\begin{align*}
    \tilde{k}_1 &=  \frac{1}{\beta}\left(\beta k_1+\dL{f}{\gamma}+\beta\dL{\lb{g}{f}}{\alpha}+\alpha\dL{g}{\gamma}-\gamma\dL{g}{\alpha}+\tilde{k}_3\left(\beta\dL{\lb{g}{f}}{\beta}+\beta\dL{g}{\gamma}-\gamma\dL{g}{\beta}\right)+\tilde{k}_2\gamma\right) \\
    &=k_1+\dL{\lb{g}{f}}{\alpha} + \frac{1}{\beta}\left(\dL{f}{\gamma}+\alpha\dL{g}{\gamma}-\gamma\dL{g}{\alpha}+\tilde{k}_3\left(\beta\dL{\lb{g}{f}}{\beta}+\beta\dL{g}{\gamma}-\gamma\dL{g}{\beta}\right)+\tilde{k}_2\gamma\right).
\end{align*}

\section{Technical lemma for the proof of \texorpdfstring{\cref{thm:characterisation_trivialisable}}{Theorem \ref*{thm:characterisation_trivialisable}}}\label{apdx:exists-diffeo}
{
The sufficiency part of the proof of \cref{thm:characterisation_trivialisable} relies on the existence of a diffeomorphism that simultaneously rectifies the distribution $\distrib{f_c,\lb{g_c}{f_c}}$ and the vector field $g_c$ as proven by the following lemma.
\begin{lemma}\label{apdx:exists-diffeo:lemma}
Consider a control-affine system $\Sigma_c=(f_c,g_c)$ given by its canonical pair and assume set $\FFF=\distrib{f_c,\lb{g_c}{f_c}}$. If, the structure functions of $\Sigma_c$ satisfy the condition \cref{eq:conditions_trivialisable} of \cref{thm:characterisation_trivialisable}, then, there exists a diffeomorphism $(x,y,w)=\phi(\xi)$ such that $\phi_*\FFF=\distrib{\vec{x},\vec{y}}$ and $\phi_*g_c=\vec{w}$. 
\end{lemma}
\begin{proof}
First, we prove that there exists smooth solutions $h$ for the system 
\begin{align*}
\dL{f_c}{h}=0,\quad \dL{\lb{g_c}{f_c}}{h}=0,\quad\textrm{and}\quad \dL{g_c}{h}=1.
\end{align*}
\noindent
We need to check three integrability conditions: 
\begin{enumerate}
\item $\dL{\lb{f_c}{\lb{g_c}{f_c}}}{h}=\dL{f_c}{\dL{\lb{g_c}{f_c}}{h}}-\dL{\lb{g_c}{f_c}}{\dL{f_c}{h}}=0$ and $\lb{f_c}{\lb{g_c}{f_c}}=0$, so $0=0$ and the first condition holds.
\item $\dL{\lb{f_c}{g_c}}{h}= \dL{f_c}{\dL{g_c}{h}}-\dL{g_c}{\dL{f_c}{h}} = 0$ and $\dL{\lb{f_c}{g_c}}{h}=0$ so $0=0$ and the second integrability condition holds.
\item $\dL{\lb{\lb{g_c}{f_c}}{g_c}}{h}= \dL{\lb{g_c}{f_c}}{\dL{g_c}{h}} - \dL{g_c}{\dL{\lb{g_c}{f_c}}{h}}=0$ and $\lb{\lb{g_c}{f_c}}{g_c}=-\varepsilon f_c-\nu\lb{g_c}{f_c}$ . Therefore $\dL{\lb{\lb{g_c}{f_c}}{g_c}}{h}=-\varepsilon\dL{f_c}{h}-\nu\dL{\lb{g_c}{f_c}}{h}=0$ and $0=0$ the third condition holds.
\end{enumerate}
Take a smooth solution $h$ of the above system, rename it $\phi_3=h$, and choose $\phi_1,\phi_2$ such that $\diff\phi_1$ and $\diff\phi_2$ annihilate $g_c$ and are independent (they exist since $g_c\neq0$). The diffeomorphism $\phi=(\phi_1,\phi_2,\phi_3)=(x,y,w)$ is such that $\phi_*\FFF=\distrib{\vec{x},\vec{y}}$ and $\phi_* g_c=\vec{w}$. 
\end{proof}
}

\end{document}